\documentclass[english,12pt]{amsart}

\usepackage{amsmath,amsfonts,bbm,amsthm}
\usepackage[utf8]{inputenc}  
\usepackage[T1]{fontenc}       
\usepackage[psamsfonts]{amssymb}

\usepackage{latexsym}
\usepackage[dvips]{epsfig}
\usepackage{dcolumn}
\usepackage{tabularx}
\usepackage{longtable}
\usepackage{graphics}
\usepackage{pstricks}
\usepackage{pst-node}
\usepackage{latexsym}
\usepackage{babel}
\usepackage{euscript}
\usepackage{dsfont}

\renewcommand{\le}{\leqslant}
\renewcommand{\ge}{\geqslant}
\newcommand{\pp}{\leqslant}
\newcommand{\pg}{\geqslant}

\newcommand{\EE}{\mathbb{E}}
\newcommand{\NN}{\mathbb{N}}
\newcommand{\PP}{\mathbb{P}}
\newcommand{\QQ}{\mathbb{Q}}
\newcommand{\RR}{\mathbb{R}}
\newcommand{\TT}{\mathbb{T}}
\newcommand{\VV}{\mathbb{V}}
\newcommand{\ZZ}{\mathbb{Z}}

\newcommand{\cD}{\mathcal{D}}
\newcommand{\cE}{\mathcal{E}}
\newcommand{\cF}{\mathcal{F}}
\newcommand{\cG}{\mathcal{G}}

\newcommand{\cI}{\mathcal{I}}

\newcommand{\cR}{\mathcal{R}}
\newcommand{\cS}{\mathcal{S}}

\newcommand{\cZ}{\mathcal{Z}}

\newcommand{\sE}{\mathsf{E}}
\newcommand{\sP}{\mathsf{P}}

\DeclareMathOperator{\id}{id}
\DeclareMathOperator{\Div}{div}

\renewcommand{\limsup}{\varlimsup}
\renewcommand{\liminf}{\varliminf}
\newcommand{\dd}[1]{\,\mathrm{d}#1}

\newcommand{\Ind}[1]{\mathbf{1}_{\left\lbrace #1 \right\rbrace}}

\newcommand{\tq}{\mathrel{/}}

\newcommand{\ens}[1]{\left\lbrace #1 \right\rbrace}

\def\Z{{\mathbb Z}}
\def\R{{\mathbb R}}
\def\w{{\omega}}
\def\indic{{{\mathbbm 1}}}
\def\demi{{1\over 2}}
\def\sERRW{$\star$-ERRW }
\def\sERRWse{$\star$-ERRW}
\def\sVRJP{$\star$-VRJP }
\def\sVRJPse{$\star$-VRJP}
	
\newtheoremstyle{bfnote}%
{}{}%
{\itshape}{}%
{\sffamily\bfseries}{. ---\hspace{0.5em}}%
{ }%
{\thmname{#1}\thmnumber{ #2}\thmnote{ (#3)}}

\theoremstyle{bfnote}
\newtheorem{defin}{Definition}[section]
\newtheorem{theorem}{Theorem}[section]
\newtheorem{corollary}[theorem]{Corollary}
\newtheorem{lemma}[theorem]{Lemma}
\newtheorem{proposition}[theorem]{Proposition}
\newtheorem{conjecture}{Conjecture}
\newtheorem{remark}{Remark}

\newenvironment{lemmabis}[1]
  {\renewcommand{\thetheorem}{\ref{#1}\textit{bis}}%
   \addtocounter{theorem}{-1}%
   \begin{lemma}}
  {\end{lemma}}

\usepackage{hyperref}

\begin{document}

\author[A. Perrel]{Adrien Perrel}
\address{Université Claude Bernard Lyon 1 \protect\\ Institut Camille Jordan, CNRS UMR 5208 \protect\\ 43, Boulevard du 11 novembre 1918, 69622 Villeurbanne cedex, France \protect\\ (\href{mailto:perrel@math.univ-lyon1.fr}{\texttt{perrel@math.univ-lyon1.fr}}).}
\author[C. Sabot]{Christophe Sabot}
\address{Université Claude Bernard Lyon 1 \protect\\ Institut Camille Jordan, CNRS UMR 5208 \protect\\ 43, Boulevard du 11 novembre 1918, 69622 Villeurbanne cedex, France \protect\\ and Institut Universitaire de France \protect\\ (\href{mailto:sabot@math.univ-lyon1.fr}{\texttt{sabot@math.univ-lyon1.fr}}).}
\title[Point of view of the particle for 2D RWDE]{Invariant measure for the process viewed from~the~particle for~2D~random~walks in~Dirichlet~environment}
\date{}

\begin{abstract}
In this paper, we consider random walks in Dirichlet random environment (RWDE) on $\Z^2$. We prove that, if the RWDE is recurrent (which is strongly conjectured when the weights are symmetric), then there does not exist any invariant measure for the process viewed from the particle which is absolutely continuous with respect to the static law of the environment. Besides, if the walk is directional transient and under condition~$\mathbf{(T')}$, we prove that there exists such an invariant probability measure if the trapping parameter verifies $\kappa>1$ or after acceleration of the process by a local function of the environment. This gives strong credit to a conjectural classification of cases of existence or non-existence of the invariant measure for two dimensional RWDE. The proof is based on a new identity, stated on general finite graphs, which is inspired by the representation of the $\star$-VRJP, a non-reversible generalization of the Vertex reinforced Jump Process, in terms of random Schr\"odinger operators. In the case of RWDE on 1D graph, the previous identity entails also a discrete analogue of the Matsumoto-Yor property for Brownian motion.
\end{abstract}

\maketitle

\section*{Introduction}


The method of the point of view of the particle, initiated by Kozlov (\cite{kozlov1985method}), has proved to be a very useful tool to investigate random walks in random conductances. In that case, the reversibility of the walk gives for free  (under a condition of integrability of the conductances) the existence of an invariant probability measure for the process viewed from the particle, which is absolutely continuous with respect to the static law. By contrast, in the case of multidimensional random walks in random environment (RWRE), where the transition probabilities are i.i.d. at each vertex, the quenched walk is not reversible and the question of the existence of such an invariant probability is a highly non-trivial question. It has been settled in three cases: in the case of general uniformly elliptic environments in the ballistic regime under Sznitman's condition~$\mathbf{(T)}$ (see \cite{peretz2022environment} where the problem is solved in all dimensions $d\ge 2$, and \cite{berger2013local} for earlier works in dimension $d\ge 4$), for balanced environments where the quenched drift is a.s. null at each vertex (\cite{Lawler82}), and in dimension $d\ge 3$ when the transition probabilities have a specific law, the Dirichlet law (see \cite{sabot2013random,bouchet2013subballistic} and \cite{sabot2017random} for an overview of the topic). In the last case, no ballisticity assumption is required and the results also include the case of walks which are not directional transient, in particular isotropic environments. For RWRE in dimension $d=2$ which are not directional transient, there is yet no available result (apart for balanced environments) and in fact there is not even a clear conjecture about existence or non-existence of such an invariant probability.

The question of the existence of an absolutely continuous invariant probability is important for several reasons. First, it gives information about trapping phenomena: indeed, if such an invariant probability exists, it implies that the process will not concentrate asymptotically on atypical regions of the environments. Besides, proving the existence of such an invariant probability may be a first step to prove diffusive behavior of the walk. Indeed, if $\pi(\w)$ is the density of the invariant distribution for environment~$\w$ with respect to the static law, then $\pi(\tau_x\w)\w(x,y)$ (where $\tau_x$ denotes the translation by~$x$) is a zero divergence field. Deep generalizations of the Kipnis-Varadhan method have been developed to understand the asymptotic behavior of random walks in divergence free fields (see e.g. \cite{komorowski2012fluctuations} for a classical reference, or \cite{kozma2017central} and \cite{toth2018quenched} for a recent and very general statements).

In this paper we consider random walks in Dirichlet environments (RWDE) on $\Z^2$. The environment is parametrized by $4$ positive parameters $(\alpha_1, \ldots, \alpha_4)$ corresponding the 4 possible steps of the walk $(e_1,e_2,e_3:=-e_1,e_4:=-e_2)$ where $(e_1,e_2)$ is the canonical base of $\Z^2$. If $d_\alpha=\sum_{i=1}^4 \alpha_i e_i$ it is known that the RWDE is directional transient if $d_\alpha\neq 0$ and that $X_n\cdot e_i$ is recurrent for $i=1,2$ if $d_\alpha=0$.  If $d_\alpha=0$ we strongly conjecture that the RWDE is recurrent, even though we are not yet able to prove it. 
The main result of the paper asserts that, under the assumption that the RWDE is recurrent, then there is no invariant measure for the process viewed from the particle which is absolutely continuous with respect to the static law. 

Let $\kappa= 2\sum_{i=1}^4 \alpha_i - \min(\alpha_1+\alpha_3, \alpha_2+\alpha_4)$ be the strength of finite traps. If $d_\alpha\neq 0$, the RWDE is directional transient and it is widely conjectured that in that case Sznitman's  condition~$\mathbf{(T')}$ is true (it is known in particular if $\vert \alpha_1-\alpha_3\vert  +\vert \alpha_2-\alpha_4\vert >1$). 
Under Sznitman's condition~$\mathbf{(T')}$, we also prove that for $\kappa >1$ there exists an absolutely continuous invariant probability for the environment viewed from the particle with $L^p$ density with respect to the static law, for all $p<\kappa$. If $\kappa\le 1$ there is no such invariant probability due to the existence of local traps where the random walk can be trapped a non-integrable amount of time. However, the effect of local traps can be canceled by an acceleration of the walk by a local function of the environment and we prove that there exists an absolutely continuous invariant probability for the accelerated process viewed from the particle. 

The proofs of the results concerning the invariant measure rely on a new identity stated on any weighted finite directed graph, which we think is of independent interest. Let us describe this identity briefly (see Section~\ref{identity} for details of notations). Consider a finite connected directed graph $\cG=(V,E)$ with positive weights $(\alpha_e)_{e\in E}$ on the edges. Take two vertices $x$~and~$y$, and assume that the weights have zero divergence outside of $x$ and $y$. Let $\w$ be a random Dirichlet environments for the weights $(\alpha_e)$ and denote by $P^\w(x\rightarrow y)$ the probability for the walk in environment $\w$ to go from $x$ to $y$ without coming back to $x$. Consider $\gamma_x$ and $\gamma_y$ two independent gamma random variables, independent of $\w$, with parameters $\sum_{\underline e=x} \alpha_e$ and $\sum_{\underline e=y} \alpha_e$. Then we prove that the law of ${\gamma_x P^\w(x\rightarrow  y) \over \gamma_y P^\w(y\rightarrow  x)}$ conditioned on ${\gamma_x P^\w(x\rightarrow  y) \gamma_y P^\w(y\rightarrow  x)}$ is explicit and given by a generalized inverse Gaussian law. 

That formula is inspired by some computations done for the Vertex Reinforced jump Process (VRJP) and more specifically for its non-reversible generalization (\cite{svrjp}), the $\star$-VRJP. In dimension 1, our formula is related to the Matsumoto-Yor formula, which concerns exponential functionals of the Brownian motion. More precisely it gives a discrete version of the Matsumoto-Yor formula.


\section{Description of the model and known results}

%

	\subsection{Random walks in Dirichlet environment on directed graphs}
	\label{rwre_def}

We first define Dirichlet environments on general graphs since it will be useful to state the main new identity. Consider $\cG=(V,E)$ a connected directed graph with bounded degree at each vertex. If $e=(x,y)$ is a directed edge, we set $\underline e:=x$ for its tail and $\overline e:=y$ for its head. The graph is said strongly connected if for any two vertices $x$ and $y$ there exists a directed path in $\cG$ from $x$ to $y$.

The set of environments on the graph is the set of transition probabilities of Markov chains on the graph:
$$ \Omega = \ens{(\omega(x, y))_{(x,y) \in E} \in (0,1]^E \tq \forall x \in V, \sum_{y \tq (x,y)\in E} \omega(x,y) = 1}. $$

%
%
For an environment $\w\in \Omega$, we denote by $\sP^\w_{x_0}$ the law of the Markov chain in the environment $\w$ starting at $x_0$, called the quenched law. More precisely, if $X_n$ is the canonical process on $V^\NN$, we have for every $(x,y)\in E^2$:
%
%
\begin{align*}
\begin{cases}
&\sP^\omega_{x_0}(X_0 = x) = \indic_{x=x_0}, \\
&\sP^\omega_{x_0}\left(X_{n+1} = y \right.\left\vert X_n=x, \{X_k, k\le n\}\right) =  \w(x,y) \indic_{(x,y)\in E},
 \quad \forall n \in \NN.
\end{cases}
\end{align*}
This probability distribution $\sP^\omega_{x_0}$ is called the \emph{quenched law} of the random walk in the environment~$\omega$.

For every subset $U \subset V$, it is convenient to define the hitting time, exit time, and return time to~$U$. They are the stopping times defined respectively by:
	\begin{align*}
	H_U &= \inf\ens{n \pg 0 \tq X_n \in U}; \\
	\bar H_U &= \inf\ens{n \pg 0 \tq X_n \not\in U}; \\
	H^+_U &= \inf\ens{n \pg 1 \tq X_{n-1} \not\in U, X_n \in U}.
	\end{align*}
We simply write $H_{x}$, $\bar H_{x}$, $H^+_{x}$, when $U$ is the singleton $U=\{x\}$.
%
%
%
%

{\bf Dirichlet distribution.} 
Given an integer $k \in \NN^*$ and $k$~positive parameters $(\alpha_1, \ldots, \alpha_k)$, the Dirichlet distribution $\cD(\alpha_1, \ldots, \alpha_k)$ is a random partition of $[0,1]$ in $k$~parts, with density:
\begin{equation}
\label{eq:Dirichlet_distrib}
\frac{\Gamma\left(\sum_{i=1}^k \alpha_i\right)}{\prod_{i=1}^{k} \Gamma(\alpha_i)} \left(\prod_{i=1}^{k} x_i^{\alpha_i-1}\right) \dd x_1 \ldots \dd x_{k-1}
\end{equation}
	restricted to the simplex $\ens{(x_i)\in(\RR_+^*)^{k} \tq \sum_{i=1}^k x_i = 1}$. It is therefore a natural generalization of the Beta distribution.
	
	Alternatively, we can construct $\cD(a_1, \ldots, a_k)$ as follows. Draw $(W_1, \ldots, W_k)$ independent random variables with respective distribution $\Gamma(a_i,1)$, and set $\Sigma = \sum_{i=1}^k W_i$. Then:
	$$\left(\frac{W_1}{\Sigma}, \ldots, \frac{W_k}{\Sigma}\right) \sim \cD(a_1, \ldots, a_k).$$

{\bf Dirichlet environment.} The Dirichlet environment on the graph $\cG
$ is parametrized by a family of positive weights $(\alpha(e))_{e\in E} \in {(\RR^*_+)}^E$ associated with the directed edges of the graph. It naturally induces weights $(\alpha(x))_{x \in V}$ for the vertices of the graph~$\cG$:
$$ \forall x \in V, \alpha(x) = \sum_{e\in E \tq \underline{e}=x} \alpha(e). $$
At each vertex $x \in V$, we draw independently a random vector $\omega_x = (\omega(e))_{\underline{e}=x}$ according to the Dirichlet law $\cD(\alpha(e), \underline{e}=x)$. The joint distribution of the $\w=(\w_x)_{x \in V}$ is denoted $\PP^{(\alpha)}$ and called the Dirichlet random environment with parameters~$\alpha$ on~$\cG$.

The annealed law is the law on $V^\NN$ given by:
$$
\sP^{(\alpha)}_{i_0}\left(\cdot\right)=\EE^{(\alpha)}\left[\sP^\w_{i_0}(\cdot)\right].
$$
When there is no risk of ambiguity, we simply write $\PP$ and $\sP_{i_0}$ for $\PP^{(\alpha)}$ and $\sP^{(\alpha)}_{i_0}$.

{\bf Divergence and time reversal.} 
The divergence is the operator $\Div: \R^E\longrightarrow \R^V$ given, for any function $(a(e))_{e\in E}$ on the edges, by:
$$ \forall x \in V, \qquad \Div(a)(x) = \sum_{e \in E \tq \underline{e}=x} a(e) - \sum_{e \in E \tq \overline{e}=x} a(e). $$

When the weights $(\alpha(e))_{e\in E}$ have zero divergence, \textit{i.~e.} $\Div(\alpha)\equiv 0$, the random Dirichlet environment with parameters $\alpha$ has a property of stability in law by time reversal. That property was key in previous works on RWDE; we recall it for sake of completeness, even though the new identity we derive in this paper does not seem to be implied by the time reversal property.

We define the reversed graph of $\cG$ as the graph $\check{\cG} = (V, \check{E})$ with same vertices but reverted edges: for every $(x, y)\in V^2$, $e = (x,y) \in E$ if and only if $\check{e}:= (y,x) \in \check{E}$. If the edges of $\cG$ are associated with positive weights $(\alpha(e))_{e\in E}$, this naturally induces positive weights $(\check\alpha(f))_{f\in \check E}$ on the edges of $\check\cG$, with $\check{\alpha}(\check{e}) = \alpha(e)$ for all $e\in E$. When $\cG$ is finite and connected, for $\w\in \Omega$ we can define the time-reversed
environment as the environment on the reversed graph $\check\cG$ defined by:
%
$$\forall (x,y) \in E, \quad \check\omega(y,x) = \frac{\pi^\omega(x)}{\pi^\omega(y)} \omega(x,y), $$
where $\pi^\w$ is the invariant measure of the environment $\w$.

\begin{lemma}[Time reversal property, \cite{sabot2011random,sabot2011reversed}]
\label{THreversal}
Assume that $\cG$ is finite and that $\Div(\alpha) \equiv 0$. Then:
$$ \omega \sim \sP^{(\alpha)} \Rightarrow \check{\omega} \sim \sP^{(\check{\alpha})}. $$
\end{lemma}

%
	
	\subsection{Dirichlet environment over $\ZZ^d$ and invariant measure for the process viewed from the particle.}
	\label{measure}
	
	
Let $(e_1, \ldots, e_d)$ be the canonical base of $\ZZ^d$ and set $e_{i+d}=-e_i$ for $1 \pp i \pp d$. We consider $\cG=(\ZZ^d,E)$ as the lattice graph, \emph{i.~e.} $(x,y)\in E$ if and only if $x-y\in \{e_1, \ldots, e_{2d}\}$. We consider weights which are invariant by translation: in other words, we fix $(\alpha_1, \ldots, \alpha_{2d})$ positive parameters, and, for $1 \pp i \pp 2d$, we consider that every edge $e = (x,y)$ such that $y-x = e_i$ is attributed a weight $\alpha(e) = \alpha_i$. It implies that for all $x\in \ZZ^d$, the transition probabilities $(\w(x,x+e_i))_{1 \pp i \pp 2d}$ are i.i.d. with distribution $\cD((\alpha_i)_{1 \pp i \pp 2d})$.

For $x\in \ZZ^d$, we denote by $\tau_x\colon \Omega\rightarrow \Omega$ the translation operator defined by
	$$ 
	\forall \w\in \Omega, \, e \in E, \qquad \tau_x\w(e)=\w(x+e).
	$$
For $\w\in \Omega$, the process $(\overline{\omega}_n)_{n \in \NN}$ of the environments viewed from the particle is defined as:
$$ \forall n \in \NN, \qquad \overline{\omega}_n = \tau_{X_n} \omega $$
where $(X_n)$ is the random walk in environment $\w$ starting at $0$, \textit{i.~e.} under law $\sP^\w_0$. The process $(\overline{\omega}_n)$ is a Markov process with spate space $\Omega$ and with generator $\cR$ given by:
	$$ \cR\phi(\omega) = \sum_{i=1}^{2d} \omega(0,e_i) \cdot \phi( \tau_{e_i} \omega), $$
for all measurable non-negative function $\phi$ over $\Omega$.


A measure $\QQ$ on $\Omega$ is said to be invariant from the point of view of the particle if it is invariant for the Markov chain~$(\overline{\omega}_n)_{n \in \NN}$, 
\emph{i.e.} if for every measurable non-negative function $f$ on $\Omega$,
$$ \int_\Omega \cR f \dd\QQ = \int_\Omega f \dd\QQ. $$
We are interested in non-trivial invariant measures $\QQ$ which are absolutely continuous with respect to the law of the environment $\PP$, \emph{i.e.} $\QQ=\rho \cdot \PP$ with $\rho$ measurable and non-negative. In that case, by ergodicity, it is known that $\rho(\w) >0$ a.s. (see e.~g. the proof of Theorem~1.2 in \cite{Sznitman_ten_lectures}) 
 and the equation above can be written in terms of $\rho$:
$$
\rho(\w)=\sum_{i=1}^{2d} \w(-e_i,0) \rho(\tau_{-e_i} \w), \;\;\; \PP\text{-a.s.}
$$
More fundamentally we are interested in absolutely continuous invariant probabilities, \emph{i.~e.} $\rho \in L^1(\PP)$ with $\int \rho \dd\PP=1$. 

	\subsection{Known results on $\ZZ^d$}
	\label{asymptotic}

We state below some results concerning transience and existence of an invariant probability. They are all consequences of the time-reversal property Lemma~\ref{THreversal}.  

\begin{theorem}[Theorem~2 of \cite{sabot2011random} and Theorem~1 of \cite{sabot2013random}]
\label{THdim3}
Consider, on $\ZZ^d$, the random walk in Dirichlet environment with parameters $(\alpha_1, \ldots, \alpha_{2d})$. Set:
$$
\kappa = \min_{1 \pp j \pp d} \left\{2 \sum_{i=1}^{2d} \alpha_i - (\alpha_j + \alpha_{j+d})\right\}.
$$
When $d \pg 3$, the RWDE is transient and more precisely,
$$\EE^{(\alpha)}\left[ (G^\w(0,0))^p\right]<\infty \hbox{ if and only if } p<\kappa,
$$
where $G^\w(x,y)=\EE_{x}\left[ \sum_{k=0}^\infty \indic_{X_k=y}\right]$ is the Green function of the quenched walk.

\noindent Besides, when $d \pg 3$,
\begin{enumerate}
\item If $\kappa > 1$, there exists a unique invariant probability measure $\QQ^{(\alpha)}$ from the point of view of the particle which is absolutely continuous with respect to~$\PP^{(\alpha)}$; moreover, $\frac{\dd\QQ^{(\alpha)}}{\dd\PP^{(\alpha)}} \in L^p$ for every $p \in [1, \kappa)$.
\item If $\kappa \pp 1$, there does not exist any invariant probability measure $\QQ^{(\alpha)}$ from the point of view of the particle which is absolutely continuous with respect to~$\PP^{(\alpha)}$.
\end{enumerate} 
\end{theorem}

The parameter $\kappa$ here introduced arises naturally when investigating the time spent by the random walk in a given edge; we discuss its signification further in section~\ref{accelerationBis}. 
%
%

Finally, by a different argument, based on the property of statistical invariance by time reversal (lemma~\ref{THreversal}), \cite{sabot2011random} and \cite{tournier2015asymptotic} have determined the regimes of directional transience for the random walk in Dirichlet environment on $\ZZ^d, d \pg 2$:
	
	\begin{theorem}[\cite{sabot2011reversed} and \cite{tournier2015asymptotic}]
	\label{THasymptotic1}
Consider, on $\ZZ^d$, the random walk in Dirichlet environment with parameters $(\alpha_1, \ldots, \alpha_{2d})$. Set:
$$
d_\alpha= \sum_{i=1}^{2d} \alpha_i e_i.
$$ 
Fix $\ell\in \RR^d\setminus\{0\}$. In any dimension $d\ge 1$, we have,
\begin{enumerate}
	\item If $d_\alpha \cdot \ell =0$, then $\sP^{(\alpha)}_0$-a.s., $\displaystyle -\infty=\liminf_{n\to\infty} X_n\cdot \ell<\limsup_{n\to\infty} X_n\cdot \ell=+\infty.$
	\item If $d_\alpha \cdot \ell >0$ (respectively $d_\alpha \cdot \ell < 0$), then $\sP^{(\alpha)}_0$-a.s., $\displaystyle \lim_{n\to\infty}X_n\cdot \ell = +\infty$ (resp. $-\infty$).
\end{enumerate}
Moreover, if  $d_\alpha \ne 0$: 
\begin{align*}
	P^{(\alpha)}_0\text{-a.s.,}\qquad \lim_{n\to\infty} \frac{X_n}{|X_n|}= \frac{d_\alpha}{|d_\alpha|}.
\end{align*}
	
	\end{theorem}

Finally, the study of the accurate behaviour of the walk in traps has given the asymptotic fluctuations of the walk around its average in the subballistic regime (\cite{poudevigne2019limit}) and in the subdiffusive regime under condition~$\mathbf{(T)}$ (\cite{perrel2024limit}).

\section{Main results}

	\subsection{A new identity for random walks in Dirichlet environment}
The following identity is the key ingredient in all the subsequent results. It is proved in section~\ref{identity}. It was inspired by some identities related to the random Schr\"odinger representation of the Vertex Reinforced Jump Process (VRJP), more precisely by its non-reversible generalization, the $\star$-VRJP (\cite{svrjp}). We give further explaination about that connection in Section~\ref{Schrodinger}.
	
\begin{lemma}[Product-ratio identity for hitting probabilities]
\label{THidentity}
Consider a finite directed and strongly connected graph $\cG = (V,E)$ and fix $i_0$ and $j_0$ two distinct vertices. Assume that the graph is endowed with positive weights $\alpha$ such that:
\begin{eqnarray}\label{cond_divergence}
\Div(\alpha) \equiv \gamma(\delta_{i_0}-\delta_{j_0}).
\end{eqnarray}
	Draw at random an environment $(\omega(e))_{e\in E}$ over the edges of $\cG$ with law $\PP^{(\alpha)}$. Draw independently from $\omega$ two independent random variables $\beta_{i_0}$ and $\beta_{j_0}$  
with Gamma laws, respectively $\Gamma(\alpha_{i_0},1)$ and $\Gamma(\alpha_{j_0},1)$.
Set
\begin{align*}
H_+ &= \beta_{i_0} \sP_{i_0}^\omega\left(H_{j_0} < H_{i_0}^+\right), &\Gamma &= \sqrt{H_+ H_-}, \\
H_- &= \beta_{j_0} \sP_{j_0}^\omega\left(H_{i_0} < H_{j_0}^+\right), &e^S &= \sqrt{\frac{H_+}{H_-}}.
\end{align*}

Then, conditionally to $\Gamma$, $S$ is a random variable with density:
$$ \frac{1}{2 K_\gamma(2\Gamma)} e^{\gamma s-2\Gamma\cosh(s)} \dd s, $$
where $K_\gamma$ is the modified Bessel function of the second kind with index~$\gamma$. Otherwise stated, conditionally on $\Gamma$, $e^S$ has the generalized inverse Gaussian law with parameters $(\gamma, 2\Gamma, 2\Gamma)$.  

In particular, if $\gamma=0$, $S$ et $-S$ follow the same law and $H_+$ et $H_-$ follow the same law.
\end{lemma}
	
	\subsection{Invariant measure viewed from the particle for random walk in Dirichlet environment on the plane}
%
	\subsubsection{Non-existence of any invariant measure when the walk is recurrent}
	
The main result of the paper concerns non-existence of an absolutely continuous invariant measure under the assumption that the RWDE is recurrent.
Remind the notations from theorem~\ref{THdim3} and~\ref{THasymptotic1}.
	
\begin{theorem}
\label{THrecurrence}
If $d_\alpha = 0$, and under the hypothesis that the random walk in Dirichlet environment~$\PP$ on $\ZZ^2$ is recurrent, 
there exists no invariant measure for the process viewed from the particle which is absolutely continuous with respect to~$\PP$.
\end{theorem}

In fact, what we really prove (and which implies Theorem~\ref{THrecurrence} above) is that there does not exist any measurable function $\w\mapsto (\pi^\w(x))_{x\in \ZZ^d}$ which is $\PP$-a.s. an invariant measure for the environment $\w$ on $\ZZ^2$ and which is stationary with respect to spacial translations.

It is strongly conjectured that the two hypotheses are redundant, \emph{i.e.} that the random walk in Dirichlet environment~$\PP$ on $\ZZ^2$ is recurrent as soon as $d_\alpha = 0$. Theorem~\ref{THrecurrence} is proved in section~\ref{recurrence}.
	
	\subsubsection{Existence of an invariant probability under condition~$\mathbf{(T')}$}
	
	We then examine the case of transient random walks in Dirichlet environment, and specifically one which satisfies condition $\mathbf{(T')}$. Before, we briefly recall Sznitman's condition $\mathbf{(T')}$.
	
Let $\hat{u}$ be a unit vector of $\RR^d$. Under the assumption that the walk is transient in the direction $\hat u$, we can define a sequence of renewal times $(T_i)_{i\in\NN^*}$
in the direction $\hat{u}$. We refer to Section~1 of \cite{sznitman1999law} for the definition, which is now very classical. We recall that $(T_i)_{i\ge 1}$ is an increasing sequence of random times which verifies that $X_n\cdot \hat u\ge X_{T_i} \cdot \hat u $ for $n\ge T_i$ and $X_n\cdot \hat u< X_{T_i} \cdot \hat u $ for $n<T_i$.

Given a unitary vector $\hat{u}$, we say that the condition~$\mathbf{(T')}$ in direction $\hat{u}$ is satisfied if:
	\begin{enumerate}
	\item $X$ is transient in the direction~$\hat u$, \textit{i.~e.}:
\begin{equation}
\label{eq:Ti} \tag{\textbf{T}i}
\langle X_n, \hat u \rangle \xrightarrow[n\to\infty]{\sP_0\text{-a.s.}} +\infty;
\end{equation}
	\item For every $a \in (0,1)$, there exists $c > 0$ such that:
\begin{equation}
\label{eq:Tii} \tag{\textbf{T}ii}
\sE_0\left[\sup_{1 \le i \le T_1} e^{c \|X_i\|^a} \right] < \infty.
\end{equation}
	\end{enumerate}

Condition~$\mathbf{(T')}$, as well as very similar conditions (named~$\mathbf{(T)}$ and~$\mathbf{(T)}_\gamma$), was introduced and first investigated in \cite{sznitman2001class, sznitman2002effective}. This bunch of conditions played since then a central role in the topic (see \textit{e.~g.}  \cite{bolthausen2002static, sznitman2003new, rassoul2009almost, yilmaz2010averaged}). 
It is conjectured that directional transience is equivalent to condition~$\mathbf{(T')}$. Partial results in this direction have been provided (\cite{berger2014effective,guerra2020proof}), but a full proof of this conjecture seems still out of reach.

In the case of Dirichlet random environments, Tournier have shown in \cite{tournier2009integrability} that condition~$\mathbf{(T)}$, hence $\mathbf{(T')}$, is satisfied as soon as
$$ \sum_{i=1}^d \vert \alpha_i - \alpha_{i+d} \vert > 1.$$ 
In \cite{Bouchet_Sabot_Ramirez}, Theorem~5, the condition~$\mathbf{(T')}$ has been proved under the condition that one of the parameter is small enough. In particular, it implies that for all $\kappa$, it is possible to find some weights $(\alpha_1, \ldots , \alpha_{2d})$ such that $\mathbf{(T')}$ is proved. Note that as explained above, it is expected that condition~$\mathbf{(T')}$ is true as soon as $d_\alpha\neq 0$.
	
\begin{theorem}
\label{THtransience}
Under the hypothesis that the random walk in Dirichlet environment~$\PP$ on $\ZZ^2$ satisfies the condition~$\mathbf{(T')}$ (and in particular that $d_\alpha \neq 0$), and that $\kappa>1$, there exists an invariant probability measure~$\QQ$ from the point of view of the particle which is absolutely continuous with respect to~$\PP$. Moreover, $\frac{\dd\QQ}{\dd\PP} \in L^p$ for every $p < \kappa$, where $\kappa$ is the parameter introduced in Theorem~\ref{THdim3}.
\end{theorem}

\begin{remark}
Note that for general environment, under the assumption of uniform ellipticity, \cite{peretz2022environment} proved the existence of such an invariant probability under the same condition~$\mathbf{(T')}$. As Dirichlet environments are not uniformly elliptic, our result is not included in that of \cite{peretz2022environment}, even though it is very comparable. Our approach is based on the identity in Lemma~\ref{THidentity} and rather different from the approach of \cite{peretz2022environment}.
\end{remark}
We prove this result in section~\ref{transience}.

\subsubsection{Accelerated random walk}
\label{acceleration}

When the walk is transient, the non-existence of an invariant probability from the point of view of the particle is linked to the effect of traps. Traps also give the interpretation of the parameter $\kappa$ introduced above (see section~\ref{accelerationBis}).

In this section, we prove that, by a sufficient acceleration of the random walk, we can counter the effect of traps, and get an absolutely continuous invariant probability from the point of view of the accelerated walk.
\begin{defin}
A function $\gamma: \Omega \longrightarrow \RR_+^*$ is called a local accelerating function if there exists a finite box $\Lambda$ of $\ZZ^d$ such that $\gamma$ is mesurable with respect to the $\sigma$-algebra generated by $\ens{\omega(x,y), (x,y) \in E_\Lambda}$. The accelerated process $(Y_t)_{t \pg 0}$ with acceleration~$\gamma$ is the continuous-time Markov chain with jump rate from $x$ to $y$:
$$ \gamma(\tau_x \omega) \omega(x,y). $$
\end{defin}

In other words, $\gamma$ is a local accelerating function if its values only depends on the environment in a finite box. The accelerated process satisfies $X_n = Y_{t_n}$, where $t_n = \sum_{k=1}^n \frac{1}{\gamma(\tau_{X_k} \omega)} E_k$, the $E_k$'s being independent exponential random variables with parameter~1: hence $Y$ is an accelerated version of $X$.


Note that we restrict ourselves to \emph{local} accelerating function (with \emph{finite} $\Lambda$). Such an accelerating function can only counter the effects of \emph{finite} traps, whose size is smaller than the size of $\Lambda$.

In dimension $d = 2$, this local effect is sufficient when the walk is transient under condition~$\mathbf{(T')}$; but it is never the case when the walk is recurrent. Namely:

	\begin{theorem}
	\label{THaccelerated}
Consider a RWDE on $\ZZ^2$. 

Then, if $d_\alpha = 0$ and if $X$ is recurrent, whatever the local accelerating function~$\gamma$, there does not exist any invariant measure from the point of view of the accelerated particle which is absolutely continuous with respect to~$\PP$.

Moreover, if $X$ satisfies condition~$\mathbf{(T')}$, then for any $p\ge 1$, there exists a local accelerating function $\gamma$ such that the accelerated process admits an absolutely continuous invariant probability $\QQ$ with respect to~$\PP$ such that $\frac{\dd\QQ}{\dd\PP} \in L^p$.
	\end{theorem}

\subsubsection{Conjectures} 
As said before, we conjecture that, in dimension~$d = 2$, the condition~$\mathbf{(T')}$ is redundant with the transience assumption, and that the walk is recurrent exactly when $d_\alpha = 0$. We can thus end this section with a global conjecture which sums up the previous remarks:

	\begin{conjecture}
	On $\ZZ^2$, consider the random walk $X$ in Dirichlet environment $\PP$ with parameters $(\alpha_1, \ldots, \alpha_4)$. Using notations from Theorems~\ref{THdim3} and~\ref{THasymptotic1}, we conjecture that:
	
	\begin{enumerate}
\item If $\kappa > 1$, there exists a unique invariant probability from the point of view of the particle if, and only if, $d_\alpha \neq 0$.

\item If $\kappa \pp 1$, there exists no invariant probability from the point of view of the particle. However, there exists a local accelerating function~$\gamma$ such that there exists a unique invariant probability from the point of view of the accelerated particle with accelerating function~$\gamma$ if, and only if, $d_\alpha \neq 0$.
	\end{enumerate}
	
	\end{conjecture}
	
\subsection{Matsumoto-Yor-type property on the line}
	
In this section, we show that in dimension 1, the product-ratio identity of Lemma~\ref{THidentity} is related to the Matsumoto-Yor (\hbox{M-Y} for short) property for exponential Brownian motion. More precisely, it gives a discrete version of M-Y property. A very similar discrete M-Y property was exhibited in \cite{arista2024matsumoto}, Theorem~1.1, with very different approach and motivation. However, although the results look extremely similar, we could not find a strict correspondence between the two. Note besides that a related discrete \hbox{M-Y} property for the drift~$\demi$ was derived in \cite{rapenne2023continuous} in relation with the VRJP. 
	
	\paragraph{\bf Matsumoto-Yor property.} We first recall \hbox{M-Y} property. Let $(W_t)_{t \pg 0}$ be a standard Brownian motion, and set, for some $\mu \in \RR$, $(E_t)_{t \pg 0}$ the drifted geometric Brownian motion:
	$$ \forall t \pg 0, \qquad E_t = e^{W_t + \mu t}. $$
	Define also the related exponential functionals $(I_t)_{t \pg 0}$ and $(Z_t)_{t \pg 0}$ as follows:
	$$ \forall t \pg 0, \qquad I_t = \int_0^t (E_t)^2 \dd t, \quad Z_t = \frac{I_t}{E_t}. $$
	We also denote by $(\cE_t)_{t \pg 0}$ and $(\cZ_t)_{t \pg 0}$ the natural filtrations adapted to $(E_t)_{t \pg 0}$ and $(Z_t)_{t \pg 0}$.
	
	\begin{theorem}[Theorem 1.6 and Proposition 1.7 of \cite{matsumoto2001analogue}]
	With the previous definitions, in the filtration $\cZ$, $Z$ is a diffusion process with known explicit infinitesimal generator. Moreover, for every $t > 0$, $\cZ_t \subsetneq \cE_t$ and the conditional distribution of $E_t$ knowing $\cZ_t$ is a generalized inverse Gaussian distribution with parameters $(\mu, 1/Z_t, 1/Z_t)$. 
	\end{theorem}
	
	\paragraph{\bf Discrete version of Matsumoto-Yor property.}
Lemma~\ref{THidentity} enables us to provide a discrete counterpart of \hbox{M-Y} property. The graph we consider is 1D segment $ \ens{0, \ldots, n}$. More precisely,
we consider the following sequence of graphs $(\cG_n)_{n \pg 1}$: for every $n \pg 1$, $\cG_n = (V_n, E_n)$ is composed of $V_n = \ens{0, \ldots, n}$ with edges between nearest neighbours $E_n = \ens{(i, i+1), (i+1,i), 0 \pp i \pp n-1}$. We endow these graphs $\cG_n$ with Dirichlet environments $\omega_n$ with parameters $(\alpha, \beta) \in (\RR_+^*)^2$ as follows: for every $i \pg 0$, we draw independently at random variables $W(i, i+1)$ and $W(i+1, i)$ with respective distributions $\Gamma(\alpha, 1)$ and $\Gamma(\beta, 1)$; we then renormalize these gamma variables, defining, for every $n \pg 0$ and $(i,j)\in E_n$:
	$$\omega_n(i,j) = \frac{W(i,j)}{\sum_{(i,k)\in E_n} W(i,k)}.$$
	(Note that, for $(i,j)=(n,n-1)$, $\omega_n(n,n-1)=1$ in $\cG_n$, but $\omega_{n+1}(n,n-1)$ has a Beta distribution in $\cG_{n+1}$.)
	
Let $(\Gamma_n)_{n \pg 1}$ and $(S_n)_{n \pg 1}$ be two sequences of random variables defined by: for every $n \pg 1$,
	\begin{align*}
	\Gamma_n &= \sqrt{W(0,1)W(n,n-1) \sP_0^{\omega_n}\left(H_n < H_0^+\right)\sP_n^{\omega_n}\left(H_0 < H_n^+\right)} \\*
	e^{S_n} &= \sqrt{\frac{W(0,1)\sP_0^{\omega_n}\left(H_n < H_0^+\right)}{W(n,n-1) \sP_n^{\omega_n}\left(H_0 < H_n^+\right)}}
	\end{align*}
	Those definitions coincide with those of Lemma~\ref{THidentity} for $\cG = \cG_n$, $i_0 = 0$ and $j_0 = n$. We can rewrite them. Define further: for every $i \pg 1$,
$$ \rho(i) = \frac{W(i,i-1)}{W(i,i+1)}. $$
Then:
	
	\begin{lemma}
	\label{THexplicit}
With the previous definition, for every $n \pg 1$:
\begin{align*}
\Gamma_n &= \sqrt{W(0,1) W(n, n-1)} \frac{\prod_{i=1}^{n-1} \sqrt{\rho(i)}}{\sum_{0 \pp k < n} \prod_{i=1}^k \rho(i)},\\
e^{S_n} &= \sqrt{\frac{W(n,n-1)}{W(0,1)}} \prod_{i=1}^{n-1} \sqrt{\rho(i)}.
\end{align*}
	\end{lemma}
	
	Let us finally also denote $(\cF_n)_{n \pg 1}$ and $(\cS_n)_{n \pg 1}$ the natural filtrations adapted to $(\Gamma_n)_{n \pg 1}$ and $(S_n)_{n \pg 1}$ respectively.
	
	\begin{theorem}[Discrete version of Matsumoto-Yor property]
	\label{THMatsYor}
	With the previous definitions, in the filtration $\cF$, $\Gamma$ is a Markov process with known explicit transition kernel. Moreover, for every $n \pg 1$, 
conditionally on $\cF_n$, $e^{S_n}$ has the inverse Gaussian law with parameters $(\alpha-\beta, 2\Gamma_n, 2\Gamma_n)$.
	\end{theorem}
	
	The previous Lemma~\ref{THexplicit} justifies the name of ``discrete Matsumoto-Yor property'' we have given to our result. Indeed, there is a formal correspondence between the exponential functional~$E_t$ of the Brownian motion of Matsumoto and Yor, and the multiplicative random walk~$S_n$ with increments $\sqrt{\rho_i}$ of our results: the integral $I_t$ of the squared functional transforms into the sum of the squared first values of $S$, and the analogous of $Z_t$ is now $\Gamma_n$. Note that the $\rho(i)$ are i.i.d. variables with beta-prime distribution $B'(\beta, \alpha)$. In particular,
$$ \EE[\ln \rho(i)] = \psi(\beta) - \psi(\alpha), \qquad \VV[\ln \rho(i)] = \psi'(\alpha) + \psi'(\beta) , $$
where $\psi$ denotes the digamma function. By taking a fine-mesh limit of our model, we could also derive a new proof of the Matsumoto-Yor property, by taking $\alpha = \frac{n}{2}$ and $\beta = \frac{n}{2} + \mu$; in this setting, the linear interpolation of $(S_{\lfloor nt \rfloor})_{t\pg 0}$ converges towards~$E$ when ${n \to \infty}$. This was done in \cite{rapenne2023continuous} for a similar discrete Matsumoto-Yor property for the drift $\mu=-1/2$. We do not prove the convergence of our discrete M-Y property to the original M-Y property since it involves some technical steps which are very similar to what has been done in \cite{rapenne2023continuous}.

\section{Proof of the product-ratio identity for hitting probabilities}
\label{identity}

%
In this section we consider the setting of Lemma~\ref{THidentity}: we consider a finite connected directed graph $\cG = (V,E)$, whose edges are endowed with positive weights $(\alpha(e))_{e\in E}$ satisfying~\eqref{cond_divergence}. We also assume, without loss of generality, that $\gamma \pg 0$.
	

	\subsection{A mixing property of Dirichlet environments}
	
	For fixed $\ens{W(e), e \in E}\in (\R^*_+)^E$, introduce random variables $(U_i)_{i \in V}$ on the vertices of $\cG$ with distribution:
	\begin{align}\label{nu_gamma}
	\nu_\gamma^{W, \cG}(\dd u) = \frac{1}{Z_\gamma(W)} e^{\gamma(u_{j_0}-u_{i_0})} \exp\left(-\sum_{(i,j)\in E} W(i,j) e^{u_j - u_i} \right) \delta_{0}(u_{i_0}) \prod_{i \in V \setminus \{i_0\}} \dd u_i
	\end{align}
 on $\ens{(u_i)_{i \in V} \in \RR^V \tq u_{i_0} = 0}$, where $Z_\gamma (W)$ is the normalizing constant.
	
	Let $(W^U(e))_{e\in E}$ be new random variables on the edges of $\cG$ defined by:
	$$ \forall (i,j)\in E, \qquad W^U(i,j) = W(i,j) e^{U_j - U_i}. $$
	That way we introduce further randomness, which enables the coupling of $W$ and $W^U$. It turns out that those two families follow the same distribution.
	
	\begin{proposition}[Mixing property]
	\label{THmixing}
	Consider a finite directed graph $\cG$ endowed with positive weights $\alpha$ such that:
	$$\Div(\alpha) \equiv \gamma(\delta_{i_0}-\delta_{j_0}).$$
and draw independent random variables $(W(e))_{e\in E}$ with respective distribution $\Gamma(\alpha(e),~1)$.

	Then $(W^U(e))_{e\in E}$ and $(W(e))_{e \in E}$ have the same distribution.
	\end{proposition}
	
	\begin{proof}
	The proof follows closely the proof of Lemma~3.7 in \cite{svrjp}.
		Let $\phi: \RR^E \longrightarrow \RR_+$ be a non-negative measurable function. We have:
\begin{align*}
\EE\left[\phi(W^U)\right] &= \int_{\RR_+^E} \int_{\RR^{V\setminus\{i_0\}}} \phi\left(\left(w_{i,j} e^{u_j-u_i} \right)\right) \frac{e^{\gamma(u_{j_0}-u_{i_0})}}{Z_\gamma(w)} 
\exp\left(-\sum_{(i,j)\in E} w_{i,j} e^{u_j - u_i}-\sum_{(i,j)\in E} w_{i,j} \right)
\\*
&\qquad\qquad\qquad\times  \prod_{(i,j)\in E}{(w_{i,j})^{\alpha(i,j)-1}dw_{i,j}\over \Gamma(\alpha(i,j))}  \prod_{i \in V\setminus\{i_0\}}\dd u_i 
\end{align*}

We proceed to the following change of variables of the $w_{i,j}$ variables, leaving $u$ unchanged:
\begin{align*}
t_{i,j} &= w_{i,j} e^{u_j - u_i}, \;\;\; (i,j)\in E. 
\end{align*}
Note that this change of variables acts on the normalization factor $Z_\gamma$ as follows: 
\begin{align*}
Z_\gamma(w)
&= Z_\gamma\left(\left(t_{i,j} e^{u_i-u_j}\right)_{(i,j)\in E}\right) \\
&= \int_{\RR_+^{V\setminus\{i_0\}}} e^{\gamma(x_{j_0}-x_{i_0})} \exp\left(-\sum_{(i,j)\in E} t_{i,j} e^{u_i-u_j} e^{x_j-x_i} \right) \prod_{i \in V\setminus\{i_0\}} \dd x_i \\
&= e^{\gamma(u_{j_0} - u_{i_0})} \int_{\RR_+^{V\setminus\{i_0\}}} e^{\gamma((x_{j_0}-u_{j_0})-(x_{i_0}-u_{i_0})} \exp\left(-\sum_{(i,j)\in E} t_{i,j} e^{(x_j-u_j)-(x_i-u_i)} \right) \prod_{i \in V\setminus\{i_0\}} \dd x_i \\
&= e^{\gamma(u_{j_0} - u_{i_0})} Z_\gamma(t).
\end{align*}
Besides, we have the following transformation by change of variable:
\begin{align*}
 \prod_{(i,j)\in E}(w_{i,j})^{\alpha(i,j)-1} \dd w_{i,j}
 =  \prod_{(i,j)\in E}e^{\alpha(i,j)(u_i-u_j)}(t_{i,j})^{\alpha(i,j)-1}\dd t_{i,j} 
 = e^{-\left<\alpha, \nabla u_{i,j}\right>}  \prod_{(i,j)\in E}(t_{i,j})^{\alpha(i,j)-1}\dd t_{i,j}.
\end{align*}
where $\nabla u_{i,j}=u_j-u_i$. But we have $-\left<\alpha, \nabla u_{i,j}\right>=\left< \Div(\alpha), u\right>= \gamma(u_{i_0}-u_{j_0})$. Hence:
$$
\prod_{(i,j)\in E}(w_{i,j})^{\alpha(i,j)-1} \dd w_{i,j}= e^{\gamma(u_{i_0}-u_{j_0})} \prod_{(i,j)\in E}(t_{i,j})^{\alpha(i,j)-1} \dd t_{i,j}.
$$

We so have:
\begin{align*}
\EE\left[\phi(W^U)\right] &= \int_{\RR_+^E} \int_{\RR^{V\setminus\{i_0\}}} \phi\left((t_{i,j})\right) \frac{e^{\gamma(u_{i_0}-u_{j_0})}}{Z_\gamma(t)} 
\exp\left(-\sum_{(i,j)\in E} t_{i,j}-\sum_{(i,j)\in E} t_{i,j}e^{u_i-u_j} \right)
\\*
&\qquad\qquad\qquad\times  \prod_{(i,j)\in E}{(t_{i,j})^{\alpha(i,j)-1}\dd t_{i,j}\over \Gamma(\alpha(i,j))}  \prod_{i \in V\setminus\{i_0\}}\dd u_i 
\end{align*}
Changing $u$ to $-u$, we get:
\begin{align*}
\EE\left[\phi(W^U)\right] &= \int_{\RR_+^E} 
\left( \int_{\RR^{V\setminus\{i_0\}}} {e^{\gamma(u_{j_0}-u_{i_0})}}
\exp\left(-\sum_{(i,j)\in E} t_{i,j}e^{u_j-u_i} \right) \prod_{i \in V\setminus\{i_0\}}\dd u_i  \right)
\\*
&\qquad\qquad\qquad\times \frac{1}{Z_\gamma(t)}  \phi\left((t_{i,j})\right) \prod_{(i,j)\in E}{e^{-t_{i,j}} (t_{i,j})^{\alpha(i,j)-1} \dd t_{i,j}\over \Gamma(\alpha(i,j))}  
\\*
&=\int_{\RR_+^E} 
   \phi\left((t_{i,j})\right) \prod_{(i,j)\in E}{(t_{i,j})^{\alpha(i,j)-1} e^{-t_{i,j}}  \dd t_{i,j}\over \Gamma(\alpha(i,j))},
\end{align*}
which concludes the proof.
\end{proof}
%
%
	
	In particular, once renormalized, $(W^U(e))_{e\in E}$ defines a new Dirichlet environment on $\cG$:
	
	\begin{corollary}[Mixing property]\label{THmixing_cor}
	Keep the notations used in proposition~\ref{THmixing}. Consider a finite directed graph $\cG$ endowed with positive weights $\alpha$ such that $\Div(\alpha) \equiv \gamma(\delta_{i_0}-\delta_{j_0})$. Let $\beta = (\beta_i)_{i \in V}$ be the normalization factors and $\omega^U = (\omega^U(e))_{e \in E}$ be the probabilities:
	$$ \forall i \in V, \quad \beta_i = \sum_{e \in E \tq \underline{e}=i} W^U(e); \qquad \forall (i,j) \in E, \quad \omega^U(i,j) = \frac{W^U(i,j)}{\beta_i}.$$

	Then, under $\PP^{(\alpha)}$, $\beta$ is a family of independent gamma variables with respective law $\Gamma(\alpha(i))$ and $\omega^U$ is a Dirichlet environment with law $\PP^{(\alpha)}$ independent from $\beta$.
	\end{corollary}
	
	In other words, $\omega$ and $\omega^U$ are coupled environment with the same law $\PP^{(\alpha)}$; but the second one has got additional structure, introduced by~$U$, which enables to reveal new properties of~$\PP^{(\alpha)}$. In the following, we will thus work with $\omega^U$ and study the variables $\Gamma^U$ and $S^U$, analogous in $\omega^U$ of $\Gamma$ and $S$ in $\omega$: it is easier to find the relation between $\Gamma^U$ and $S^U$ than between $\Gamma$ and $S$ precisely because this relation lies in the additional structure of $\omega^U$.
	
	\subsection{Proof of the product-ratio identity}

	Since we know the joint distribution of $U$, we can deduce the joint distribution of $(U_{i_0}, U_{j_0}, (\beta_i)_{i \in \tilde V})$, where we set
	$$
	\tilde V=V\setminus \{i_0,j_0\}.
	$$
	
	For that, the following notations for matrices will prove to be useful. Define the matrices:
$$ B = \text{diag}(\beta_i)_{i \in V}, \quad e^U=(e^{U_i})_{i\in V}, 
\quad W = (W(i,j))_{(i,j)\in V^2}, \quad H_\beta = B-W.
$$
with the convention that $W(i,j)=0$ if $(i,j)\notin E$.

Given $M = (m_{i,j})_{(i,j)\in I\times J}$ any matrix and $I_1 \subset I, J_1 \subset J$, we denote the restriction of $M$ to $I_1 \times J_1$ the following way:
$$ M_{I_1, J_1} = (m_{i,j})_{(i,j) \in I_1 \times J_1}.$$
In particular, 
we define short-hand notations for restrictions to~$\tilde V$:
\begin{eqnarray}\label{notations_check} 
\hat B = B_{\tilde V, \tilde V}, 
\quad \hat W = W_{\tilde V, \tilde V}, \quad \hat H_\beta = (H_\beta)_{\tilde V, \tilde V} 
\end{eqnarray}
We write $\hat H_\beta >0$ to mean that all the eigenvalues of $\hat H_\beta$ have positive real parts. In that case $\hat H_\beta$ is a $M$-matrix since its off-diagonal coefficients are non-positive (see \cite{berman1994matrices} Chapter~6, or appendix~A in \cite{svrjp}). When $\hat H_\beta >0$, we can define
$$
\hat G_\beta= (\hat H_\beta)^{-1},
$$
which has positive coefficients (see \cite{berman1994matrices} Chapter 6 Theorem~2.3~$(N_{38})$, or appendix~A in \cite{svrjp}).

\begin{lemma}
\label{THdensity}
Conditionally to $W$, $(U_{i_0}, U_{j_0}, (\beta_i)_{i \in \tilde V})$ is a random vector with distribution on $\{0\} \times \RR \times \RR^{\tilde V}$:
\begin{align*}
&\frac{\Ind{\hat{H}_b > 0}  }{Z_\gamma (W) \left\vert \det \hat{H}_b \right\vert} \exp\left(\gamma(u_{j_0} - u_{i_0}) - \sum_{i\in\tilde V} b_i\right) \\*
&\qquad \exp\left(- \check{W}_+ e^{u_{j_0}} - \check{W}_- e^{-u_{j_0}} - \check{W}_{++} - \check{W}_{--} \right)    \delta_{0}(u_{i_0}) \dd u_{j_0}\left(\prod_{i\in\tilde V} \dd b_i\right),
\end{align*}
where:
\begin{align*}
\check{W}_{+} &= W_{i_0,\tilde V} \hat G_b W_{\tilde V, j_0} + W(i_0, j_0), & \check{W}_{++} &= W_{i_0,\tilde V} \hat G_b W_{\tilde V, i_0}, \\
\check{W}_{-} &= W_{j_0,\tilde V} \hat G_b W_{\tilde V, i_0} + W(j_0, i_0), & \check{W}_{--} &= W_{j_0,\tilde V} \hat G_b W_{\tilde V, j_0}.
\end{align*}
where we recall that $W(i_0, j_0) = 0$, resp. $W(j_0, i_0)=0$, if $(i_0,j_0)\not\in E$, resp. $(j_0,i_0)\not\in E$.
\end{lemma}
	
\begin{proof}
We must proceed to the change of variables $(U_i)_{i \in V\setminus\{i_0\}}$ to $(U_{j_0}, (\beta_i)_{i \in \tilde V})$. Remind that
for all $i\in \tilde V$, $ \beta_i = \sum_{(i,j)\in E} W(i,j) e^{U_j - U_i},$
and that conditionally on $W$,  $U$ has distribution $\nu_\gamma^{\cG, W}$ given by \eqref{nu_gamma}. 
We need to make the change of variable: 
\begin{eqnarray}\label{change_var}
\R^{V\setminus \{i_0\}} &\to & \R\times \{ (b_i)_{i\in \tilde V}, \hat H_b >0\}
\\
\nonumber (u_i)_{i\in V\setminus \{0\}} &\mapsto &  (u_{j_0}, (b_i)_{i\in \tilde V})
\end{eqnarray}
where we have defined, considering that $u_{i_0}=0$,
\begin{eqnarray}\label{b_u}
b_i=\sum_{(i,j)\in E} W(i,j) e^{u_j - u_i}.
\end{eqnarray}

Let us first prove that the change of variable is bijective.
From \eqref{b_u} we get:
$$ \left(H_b e^u\right)_{\tilde V} = 0. $$
We isolate the contributions of $\tilde V$ and $\{i_0, j_0\}$ in the previous relation. For $i \in \tilde V$, we have:
\begin{align*}
0 = [H_b e^u]_{i,1}
&= \sum_{k \in V} [H_b]_{i,k} e^{u_k} \\*
&= \sum_{k \in \tilde V} [H_b]_{i,k} e^{u_k}
+ [H_b]_{i, i_0} e^{u_{i_0}} + [H_b]_{i, j_0} e^{u_{j_0}} \\*
&= [\hat H_b (e^u)_{\tilde V}]_{i,1} - W(i,i_0) - W(i,j_0).
\end{align*}
In matrix form, we get:
$$ \hat H_b ((e^u)_{\tilde V}) = W_{\tilde V, i_0} + W_{\tilde V, j_0} e^{u_{j_0}}. $$
Since $e^u$ is a vector with positive coefficients it implies that $\hat H_b$ is a $M$-matrix, i.e. $H_b >0$ (see \cite{berman1994matrices} or appendix~A in \cite{svrjp}). 
The inverted relation gives:
\begin{align}\label{b_to_u}
 (e^u)_{\tilde V} = \hat G_b \left(W_{\tilde V, i_0} + W_{\tilde V, j_0} e^{u_{j_0}}\right).
\end{align}
Reciprocally, if $(b_i)_{i\in \tilde V}$ satisfies $\hat H_b>0$, then $\hat G_b$ has positive coefficients and we can define $(e^u)_{\tilde V}$ by the previous formula \eqref{b_to_u}. This proves that the change of variables \eqref{change_var} is bijective and it is obviously $C^1$.

From \eqref{b_to_u} we also deduce that,
$$
\sum_{j \tq (i_0,j)\in E} W(i_0,j) e^{u_j}= W(i_0,j_0) e^{u_{j_0}} + W_{i_0,\tilde V} (e^u)_{\tilde V}= \check W_{++} + \check W_{+} e^{u_{j_0}},
$$
and similarly,
$$
\sum_{j \tq (j_0,j)\in E} W(j_0,j) e^{u_j-u_{j_0}}= \check W_{--} + \check W_- e^{-u_{j_0}}.
$$

Coming back to the density of $U$, we write:
\begin{align*}
&\frac{1}{Z_\gamma(W)} e^{\gamma(u_{j_0}-u_{i_0})} \exp\left(-\sum_{(i,j)\in E} W(i,j) e^{u_j - u_i} \right) \\*
&\quad = \frac{1}{Z_\gamma(W)} e^{\gamma(u_{j_0}-u_{i_0})} \exp\left(-\sum_{i\in\tilde V} b_i - \sum_{j \tq (i_0,j)\in E} W(i_0,j) e^{u_j} - \sum_{j \tq (j_0,j)\in E} W(j_0,j) e^{u_j-u_{j_0}} \right) \\
&\quad = \frac{1}{Z_\gamma(W)} e^{\gamma(u_{j_0}-u_{i_0})} \exp\left(-\sum_{i\in\tilde V} b_i - \check W_+ e^{u_{j_0}} - \check W_- e^{-u_{j_0}} - \check W_{++} - \check W_{--} \right).
\end{align*}

All that remains is to compute the Jacobian determinant of this change of variable. 
We differentiate equation \eqref{b_u}, which gives, since $u_{i_0}=0$,
\begin{align*}
& \forall i \in \tilde V, \dd b_i = \sum_{j\neq i_0, (i,j)\in E} W(i,j) e^{u_j - u_i} (\!\dd u_j - \dd u_i) \\
\Leftrightarrow &\dd b = 
\left( \text{diag}(e^{u})^{-1} \hat H_b \; \text{diag}(e^{u})\right)_{\tilde V, V\setminus{i_0}} \dd u
\end{align*}
where $\text{diag} (e^{u})$ is the diagonal matrix with diagonal coefficients $(e^{u_i})_{i\in V}$.
Finally:
$$ \dd u_{j_0} \prod_{i \in \tilde V} \dd b_i = \left\vert \det \hat{H}_\beta \right\vert \prod_{i \in V\setminus \{i_0\} } \dd u_i. $$

Using all the previous computations, we finally get: for all $\phi: \RR \times \RR_+^{\tilde V} \longrightarrow \RR_+$ non-negative measurable function,
\begin{align*}
&\EE\left[\phi\left(U_{j_0}, (\beta_i)_{i \in\tilde V}\right)\right] \\*
= & \int_{\R\times \ens{\hat H_b > 0}} \phi\left(u_{j_0}, (b_i)_{i \in\tilde V}\right) 
\exp\left(\gamma(u_{j_0}-u_{i_0}) - \sum_{i\in\tilde V} b_i - \check W_+ e^{u_{j_0}} - \check W_- e^{-u_{j_0}} - \check W_{++} - \check W_{--}\right)\\*
&\quad \quad \quad \quad \times\frac{1}{Z_{\gamma}(W) \left\vert \det \hat{H}_b \right\vert}  \left( \prod_{i\in\tilde V} \dd b_i \right) \dd u_{j_0}.
\end{align*}
\end{proof}

In the previous lemma, the role of $U_{j_0}$ and $-U_{j_0}$ is not symmetric. We symmetrize it.

\begin{lemma}
\label{THlaw}
Using notations from lemma~\ref{THdensity}, let:
$$ \Gamma^U = \sqrt{\check{W}_+ \check{W}_-}, \qquad S^U = U_{j_0} - \ln \frac{\check{W}_-}{\check{W}_+}. $$

Conditionally to $\Gamma^U$, $S^U$ is a random variable with density on~$\RR$:
$$ \frac{1}{2 K_\gamma(2\Gamma^U)} e^{\gamma s -2\Gamma^U\cosh(s)} \dd s. $$
\end{lemma}

\begin{proof}
We deduce from the previous lemma that the density of $U_{j_0}$ conditionally to $W$ and $(\beta_i)_{i \in \tilde V}$ is proportional to:
$$ \exp\left(\gamma u_{j_0} - \check{W}_+ e^{u_{j_0}} - \check{W}_- e^{-u_{j_0}}\right) \dd u_{j_0}. $$

By the change of variable $s=u_{j_0}- \ln \frac{\check{W}_-}{\check{W}_+}$, the density of~$S^U$ conditionally to $W$ and $(\beta_i)_{i \in \tilde V}$ is proportional to:
$$ \exp\left(\gamma s - \Gamma^U (e^s + e^{-s})\right) \dd s. $$
We recognize the integrand of a modified Bessel function of the second kind:
$$ \forall \alpha \in \RR, \forall x > 0, \qquad K_\alpha(x) = \int_0^\infty e^{-x \cosh(t)} \cosh(\alpha t) \dd t = \frac{1}{2} \int_{-\infty}^\infty e^{\alpha t} e^{-x \cosh(t)} \dd t. $$
Therefore, conditionally to $W$ and $(\beta_i)_{i \in \tilde V}$, $S^U$ is a continuous random variable with density:
$$ \frac{1}{2 K_\gamma(2\Gamma^U)} e^{\gamma s -2\Gamma^U\cosh(s)} \dd s. $$

Note that, at fixed~$s$, this density is a measurable function of~$\Gamma^U$, which itself is a measurable function of $W$ and $(\beta_i)_{i \in \tilde V}$. The previous density is thus also the density of $S^U$ given $\Gamma^U$.
\end{proof}
	
	\subsection{Link with the invariant measure}

In order to show that Lemma~\ref{THlaw} implies Lemma~\ref{THidentity}, we interpret now $\Gamma^U$ and $S^U$ in terms of random walks.

\begin{lemma}
\label{THinterpretation}
Using notations from lemma~\ref{THlaw}:
\begin{align*}
\check{W}_+ &= e^{-U_{j_0}} H^U_+ &\text{where } H_+^U &= \beta_{i_0} \sP_{i_0}^{\omega^U} \left(H_{j_0} < H_{i_0}^+\right), \\*
\check{W}_- &= e^{U_{j_0}} H^U_- &\text{where } H_-^U &= \beta_{j_0} \sP_{j_0}^{\omega^U} \left(H_{i_0} < H_{j_0}^+\right), \\*
\Gamma^U &= \sqrt{H^U_+ H^U_-}, \\*
e^{S^U} &= \sqrt{\frac{H^U_+}{H^U_-}} = \sqrt{\frac{\beta_{i_0} \pi^{\omega^U}(j_0)}{\beta_{j_0} \pi^{\omega^U}(i_0)}}, 
\end{align*}
where $\pi^{\w^U}$ is the invariant probability for the environment $\w^U$.
\end{lemma}

\begin{proof} We begin by finding a link between $\hat G_\beta$ and the random walk in environment~$\omega^U$.

We define same notations as in \eqref{notations_check}, replacing $W(e)$ by $W^U(e)$:
$$ \hat W^U = (W^U(i,j))_{(i,j)\in \tilde V} = \left(W(i,j) e^{U_j-U_i}\right)_{(i,j)\in \tilde V^2}, \quad \hat H^U_\beta = \hat B - \hat W^U, \qquad \hat G^U_\beta = (\hat H^U_\beta)^{-1}. $$
Those matrices can be rewritten using the transition matrix of the random walk in environment~$\omega^U$. Indeed, if we denote this transition matrix and its restriction as follows:
$$ \Omega^U = (\omega^U(i,j))_{(i,j)\in V^2} \qquad \text{and} \qquad \hat\Omega^U = (\omega^U(i,j))_{(i,j)\in \tilde{V}^2}, $$
we have:
$$ \Omega^U = B^{-1} {W^U} \qquad \text{and} \qquad \hat G^U_\beta = (\id_{\tilde V, \tilde V} - \hat\Omega^U)^{-1} \hat B^{-1}. $$
The matrix $(\id_{\tilde V, \tilde V} - \hat\Omega^U)^{-1}:= \hat G_{\omega^U} $ is the Green function of the Markov chain in environment $\omega^U$ killed at its exit of $\tilde V$, and we have:
$$ \hat G^U_\beta = \hat G_{\omega^U} B^{-1}. $$
Finally, for $i, j \in \tilde V$,
$$ \hat G_\beta (i,j)= e^{U_i} (\hat G_{\omega^U} B^{-1})_{i,j} e^{-U_j}. $$

We can now examine $\check{W}_+$. Thanks to the last identity, we get:
\begin{align*}
\check{W}_+ &= W_{i_0, \tilde V} \hat G_\beta W_{\tilde V, j_0} + W(i_0,j_0) \\
&= \sum_{\substack{k\neq j_0 \tq (i_0, k) \in E \\ l \neq i_0 \tq (l, j_0) \in E}} W(i_0, k) \frac{e^{U_k-U_l}}{\beta_l}  \hat G_{\omega^U}(k,l) W(l, j_0) + W(i_0,j_0).
\end{align*}
We express this last sum using only $\omega^U$:
\begin{align*}
&\check{W}_+ 
= \sum_{\substack{k\neq j_0 \tq (i_0, k) \in E \\ l \neq i_0 \tq (l, j_0) \in E}} \beta_{i_0} e^{U_{i_0}} \frac{W(i_0, k) e^{U_k - U_{i_0}}}{\beta_{i_0}} \hat G_{\omega^U}(k,l) \frac{W(l, j_0) e^{U_{j_0}-U_l}}{\beta_l} e^{-U_{j_0}} + W(i_0,j_0) \\
&= \beta_{i_0} e^{-U_{j_0}} \left[ \sum_{\substack{k\neq j_0 \tq (i_0, k) \in E \\ l \neq i_0 \tq (l, j_0) \in E \\ n\in\NN}} \sP_{i_0}^{\omega^U}\left(X_1 = k, (X_2, \ldots, X_n) \in \tilde V^{n-1}, X_{n+1}=l, X_{n+2}=j_0\right) + \omega^U(i_0,j_0) \right] \\
&= \beta_{i_0} e^{-U_{j_0}} \sP_{i_0}^{\omega^U} \left(H_{j_0} < H_{i_0}^+\right).
\end{align*}
Hence, we get:
$$ \check{W}_+ = e^{-U_{j_0}} H^U_+. $$

The proof is the same as for $\check{W}^U_-$.
\end{proof}

The previous lemmas prove Lemma~\ref{THidentity}. Indeed, Lemma~\ref{THlaw} gives the conditional law of $S^U$ knowing $\Gamma^U$, and Lemma~\ref{THinterpretation} shows that $\Gamma^U$ and $S^U$ are in fact constructed from $\omega^U$ the same way as $\Gamma$ and $S$ are constructed from $\omega$. Moreover, we proved in Lemma~\ref{THmixing} that $\omega$ and $\omega^U$ have the same law. Therefore, Lemma~\ref{THlaw} keeps being valid replacing $S^U$ and $\Gamma^U$ by $S$ and $\Gamma$. This is exactly Lemma~\ref{THidentity}.

	\subsection{Relation with the random Schr\"odinger representation of the $\star$-VRJP}\label{Schrodinger}
The proof of the product-ratio formula for $\gamma=0$ is inspired by some computations done in \cite{bacallado2023edge,svrjp} for the $\star$-Edge Reinforced random Walk (\sERRWse) and the $\star$-Vertex Reinforced Jump Process (\sVRJPse). The \sERRW and the \sVRJP are non-reversible generalizations of the now classical ERRW and VRJP. They are discrete time, resp. continuous time, self-interacting processes living on 
the vertex set $\tilde V$ of a directed graph $\tilde \cG=(\tilde V, \tilde E)$ endowed  with an involution $\star: \tilde V\to \tilde V$ such that
$$
(i,j)\in \tilde E \hbox{ if and only if } (j^\star, i^\star)\in \tilde E.
$$
We refer to \cite{svrjp} for the definitions of the \sERRW and the \sVRJPse. Consider now a RWDE on a graph $\cG=(V,E)$ with parameters $(\alpha(e))$. We known that the annealed law of the RWDE is a directed Edge Reinforced Random Walk. Assume that $\hbox{div}(\alpha)=0$. In that case, the RWDE can be considered in a non-trivial way as a \sERRWse. Indeed, consider $\check{\cG}=(\check V\simeq V, \check E)$ obtained by reversing each edge of $E$. Then define $\tilde \cG$ as the union of the two graphs $\cG$ and $\check{\cG}$ glued at a single vertex $i_0$. The $\star$-involution is just defined as the application which sends any vertex in $V$ to its copy in $\check V$. Hence we have $i_0=i_0^\star$ by the gluing of the two graphs.

Consider now the \sERRW on $\tilde \cG$ starting at $i_0$ up to a return time to $i_0$. It makes excursions from $i_0$ in $\cG$ and in $\check{\cG}$. If we map all the excursion in $\check{\cG}$ to reversed excursions in $\cG$, then the resulting path has the same distribution as the annealed law of the RWDE since $\hbox{div}(\alpha)=0$, by the time-reversal property (see e.g. formula~(3.1) in \cite{sabot2011reversed}).

The proof of Lemma~\ref{THidentity} is inspired by Lemma~3.7, Theorem~4.2 and Proposition~6.7 of \cite{svrjp} applied to the previous graph. Indeed, if $W$ is the $\star$-symmetric conductances of the \sVRJP on the graph $\tilde \cG$, then the mixing field of the \sVRJP takes a simpler form: in this case $(U_i)_{i\in \tilde V}$ is such that $U_{|V}$ has distribution proportional to
\begin{align}\label{u_star}
e^{-\sum_{(i,j) \in E} W_{i,j} e^{u_j-u_i}},
\end{align}
with $u_{i_0}=0$, and $U_{|\check V}$ is such that, 
$$
\forall i\in V, \;\;\; \sum_{j\in V, (i,j)\in E} W_{j,i} e^{U_{j^\star}-U_{i^\star}}= \sum_{j\in V, (i,j)\in E} W_{i,j} e^{U_j-U_i}.
$$
Remark that the distribution \eqref{u_star} coincides with the distribution we introduced in \eqref{nu_gamma}. In fact, with our special choice of graph, then $U_{| V}$ has also the same distribution as $A_{|V}$ where $A$ is the initial random local time which makes the \sVRJP exchangeable (but $U$ and $A$ have not the same law since $A^\star=-A$ which is not the case for $U$). Hence, the proof of our Proposition~\ref{THmixing} is essentially the same as the proof of lemma~3.7 in \cite{svrjp}. 

In Part~II of \cite{svrjp} a random potential $(\beta_i)_{i\in \tilde V}$ was introduced. It can roughly (in fact except at $i_0$) be related to the mixing field $U$ by
$$
\beta_i=\sum_{j\in V, (i,j)\in E} W_{i,j} e^{U_{j}-U_{i}},
$$
and hence $\beta_{i^*}=\beta_i$. It is proved in Proposition~6.7 of \cite{svrjp} for any subset $U=U^\star\subset \tilde V$, that the law of $\beta_{U}$ conditioned on $\beta_{U^c}$ can be computed. Our Lemma~\ref{THdensity} is inspired by that identity in the case $U=\{i_0,j_0\}$. Finally, the change of variables in Lemma~\ref{THlaw} can also be interpreted in the sense that the random variables $S^U$ in Lemma~\ref{THlaw} coincides with the random variable $\demi(U_{j_0}+U_{j_0^\star})$.

\section{Proof of discrete Matsumoto-Yor property on the line}
\label{MatsYor}

	This section is devoted to the proof of Theorem~\ref{THMatsYor}, our discrete version of Matsumoto-Yor property. The following proof relies heavily on the proof of the product-ratio identity, what prompts us into giving this proof now. Recall the frame defined before Theorem~\ref{THMatsYor}.
	
	\subsection{Closed expressions of $\Gamma_n$ and $S_n$}
	
	We first prove the expressions of $\Gamma$ and $S$ given in Lemma~\ref{THexplicit}.
	
	\begin{proof}
	Let us recall the result of gambler's ruin (see for example equation~2.1.4 in \cite{zeitouni2004random}): for every $0 \pp a \pp x \pp b$,
	\begin{align*}
	\sP_x^{\omega_n}\left(H_b < H_a\right) &= \frac{\sum_{a \pp k < x} \prod_{i=1}^{k} \rho(i)}{\sum_{a \pp k < b} \prod_{i=1}^{k} \rho(i)}, \\
	\sP_x^{\omega_n}\left(H_a < H_b\right) &= \frac{\sum_{x \pp k < b} \prod_{i=1}^{k} \rho(i)}{\sum_{a \pp k < b} \prod_{i=1}^{k} \rho(i)}.
	\end{align*}
In particular, we have
	\begin{align*}
	\sP_0^{\omega_n}\left(H_n < H_0^+\right) &= \omega_n(0, 1) \sP_{1}^{\omega_n}\left(H_n < H_0\right) =  \frac{W(0,1)}{\beta_n(0)} \frac{1}{\sum_{0 \pp k < n} \prod_{i=1}^{n} \rho(i)}, \\
	\sP_n^{\omega_n}\left(H_0 < H_n^+\right) &= \omega_n(n, n-1) \sP_{n-1}^{\omega_n}\left(H_0 < H_n\right) = \frac{W(n, n-1)}{\beta_n(n)} \frac{\prod_{i=1}^{n-1} \rho(i)}{\sum_{0 \pp k < n} \sum_{i=1}^{k} \rho(i)}.
	\end{align*}
	Injecting these expressions in the definitions of $\Gamma_n$ and $S_n$ we conclude:
\begin{align*}
\Gamma_n &= \sqrt{W(0,1) W(n, n-1)} \frac{\prod_{i=1}^{n-1} \sqrt{\rho(i)}}{\sum_{0 \pp k < n} \prod_{i=1}^k \rho(i)},\\
e^{S_n} &= \sqrt{\frac{W(n,n-1)}{W(0,1)}} \prod_{i=1}^{n-1} \sqrt{\rho(i)}.
\end{align*}
	\end{proof}
	
	\begin{lemma}
	\label{THrecursion}
	The sequence $(\Gamma_n)_{n \pg 1}$ is defined by the recursion:
	\begin{align*}
	\begin{cases}
	\Gamma_{n+1} &= \Gamma_n \dfrac{\sqrt{W(n,n+1)W(n+1,n)}}{W(n,n+1) + \Gamma_n e^{-S_n}} \quad \text{if } n \pg 1 \\
	\Gamma_1 &= \sqrt{W(0,1)W(1,0)}
	\end{cases}
	\end{align*}
	\end{lemma}
	
	\begin{proof}
	The following trivial relations will prove to be useful: for every $n \pg 1$ and $k \pg 0$,
	\begin{align*}
	\sP_0^{\omega_n}\left(H_n < H_0^+\right) &= \sP_0^{\omega_{n+k}}\left(H_n < H_0^+\right), \\
	\sP_n^{\omega_n}\left(H_0 < H_n^+\right) &= \sP_{n-1}^{\omega_{n+k}}\left(H_0 < H_n\right).
	\end{align*}
	
	By definition, for every $n \pg 1$:
	$$ \Gamma_{n+1} = \sqrt{\beta_{n+1}(0) \beta_{n+1}(n+1) \sP_0^{\omega_{n+1}}\left(H_{n+1} < H_0^+\right) \sP_{n+1}^{\omega_{n+1}}\left(H_0 < H_{n+1}^+\right)}.$$
	We need to work out $\sP_0^{\omega_{n+1}}\left(H_{n+1} < H_0^+\right)$ and $\sP_{n+1}^{\omega_{n+1}}\left(H_0 < H_{n+1}^+\right)$ in order to get a recursion relation.
	
	For the first one, we have:
	\begin{align*}
\sP_0^{\omega_{n+1}}\left(H_{n+1} < H_0^+\right) 
&= \sP_0^{\omega_{n+1}}\left(H_{n} < H_0^+\right) \sP_n^{\omega_{n+1}}\left(H_{n+1} < H_0\right) \\
&= \sP_0^{\omega_n}\left(H_{n} < H_0^+\right) \sP_n^{\omega_{n+1}}\left(H_{n+1} < H_0\right) \\
&= \sP_0^{\omega_n}\left(H_{n} < H_0^+\right) \frac{\omega(n,n+1)}{1-\omega(n,n-1) \sP_{n-1}^{\omega_{n+1}}\left(H_n < H_0\right)} \\
&= \sP_0^{\omega_n}\left(H_{n} < H_0^+\right) \frac{\omega(n,n+1)}{1-\omega(n,n-1) \sP_n^{\omega_n}\left(H_n^+ < H_0\right)} \\
&= \frac{\omega(n,n+1) \sP_0^{\omega_n}\left(H_{n} < H_0^+\right)}{1-\omega(n,n-1) \left[1-\sP_n^{\omega_n}\left(H_0 < H_n^+\right)\right]}.
	\end{align*}
	
	For the second one, we have:
	\begin{align*}
&\sP_{n+1}^{\omega_{n+1}}\left(H_0 < H_{n+1}^+\right)
= \sP_n^{\omega_{n+1}}\left(H_0 < H_{n+1}\right) \\
&\qquad= \omega(n, n-1) \sP_{n-1}^{\omega_{n+1}}\left(H_0 < H_{n+1}\right) \\
&\qquad= \omega(n, n-1) \left[\sP_{n-1}^{\omega_{n+1}}\left(H_0 < H_{n+1}, H_n < H_0\right) + \sP_{n-1}^{\omega_{n+1}}\left(H_0 < H_{n+1}, H_n > H_0\right)\right] \\
&\qquad= \omega(n, n-1) \left[\sP_{n-1}^{\omega_{n+1}}\left(H_n < H_0\right)\sP_n^{\omega_{n+1}}\left(H_0 < H_{n+1}\right) + \sP_{n-1}^{\omega_{n+1}}\left(H_0 < H_n\right)\right] \\
&\qquad= \omega(n, n-1) \left[\sP_n^{\omega_n}\left(H_n^+ < H_0\right)\sP_{n+1}^{\omega_{n+1}}\left(H_0 < H_{n+1}^+\right) + \sP_n^{\omega_n}\left(H_0 < H_n^+\right)\right],
	\end{align*}
hence:
	\begin{align*}
\sP_{n+1}^{\omega_{n+1}}\left(H_0 < H_{n+1}^+\right)
&= \frac{\omega(n, n-1)\sP_n^{\omega_n}\left(H_0 < H_n^+\right)}{1-\omega(n, n-1)\sP_n^{\omega_n}\left(H_n^+ < H_0\right)} \\
&= \frac{\omega(n, n-1)\sP_n^{\omega_n}\left(H_0 < H_n^+\right)}{1-\omega(n, n-1)\left[1-\sP_n^{\omega_n}\left(H_0 < H_n^+\right)\right]}.
	\end{align*}
	
	Returning to $\Gamma_{n+1}$, we get:
	\begin{align*}
\Gamma_{n+1}
&= \sqrt{W(0,1)W(n+1,n) \frac{\omega(n,n+1)\omega(n,n-1)\sP_0^{\omega_n}\left(H_n < H_0^+\right)\sP_n^{\omega_n}\left(H_0 < H_n^+\right)}{\left[1-\omega(n,n-1)\left(1-\sP_n^{\omega_n}\left(H_0 < H_n^+\right)\right)\right]^2}} \\
&= \Gamma_n \times \frac{\sqrt{W(n,n+1)W(n,n-1)}}{W(n,n+1)-W(n,n-1)P_n^{\omega_n}\left(H_0 < H_n^+\right)}.
	\end{align*}
	
	Finally, note that, for every $n \pg 1$:
	$$ \Gamma_n e^{-S_n} = W(n,n-1) \sP_n^{\omega_n}\left(H_0 < H_n^+\right), \qquad \Gamma_n e^{S_n} = W(0,1) \sP_0^{\omega_n}\left(H_n < H_0^+\right). $$
	
	This gives directly the recursion relation. The initialization is trivial.
	\end{proof}
	
	Conditionally to $\Gamma_n$, note that the distribution $S_n$ is known and only depends on $\Gamma_n$ (by Lemma~\ref{THidentity}), and that $W(n,n+1)$ and $W(n+1,n)$ are independent on~$\Gamma_n$ (the latter variable being measurable with respect to $\omega_n$, which is it independent on the former variables). We have so determined the distribution of $\Gamma_{n+1}$ conditionally to $\Gamma_n$.
	
	\subsection{Markovian character of $\Gamma_n$}
	
	However, to prove the Markovian character of $(\Gamma_n)$, we need to ensure that the conditional distribution of $\Gamma_{n+1}$ (or, equivalently, of $S_n$) knowing $\Gamma_n$ is also the conditional distribution of $\Gamma_{n+1}$ (respectively $S_n$) knowing $\cF_n = \sigma\left(\Gamma_1, \ldots, \Gamma_n\right)$. We naturally return to the proof of the product-ratio identity to check it: please refer to section~\ref{identity} for the definitions of the variables used in that proof.
	
As we have seen in the proof of the product-ratio formula, we have computed the law of $S$ conditioned to the sigma field $\sigma(W, (\beta_{i})_{i\notin \{i_0,  j_0\}})$, which is larger than the sigma field of $\Gamma$. The Markovian character of $\Gamma_n$ will come from that property. The idea is to replace the random variables $W$ by the variables $W^U$, which have same law but with an extra structure. One difficulty is that the sequence $(\Gamma_n)$ is defined for the sequence of graphs $(\cG_n)$ and that the variable $U$ is defined for each graph $\cG_n$. The first step is to prove that the random variables $U$ are compatible and can thus be coupled for all graphs $\cG_n$.

	
%

\begin{lemma}
\label{THcompatible}
The variables $U$ defined for Proposition~\ref{THmixing} are compatible on the half-line, in the sense that we can define a sequence of random variables $(U_i)_{i\in \NN}$ such that, for all $n \pg 1$, $(U_0, \ldots, U_{n})$ has law $\nu_\gamma^{W, \cG_{n}}$.
\end{lemma}
	
%
%
	
	\begin{proof}
The distribution $\nu_\gamma^{W, \cG_{n}}$ can be written under the form, with $u_0=0$:
$$
\nu_\gamma^{W, \cG_{n}}(du)= \exp\left(-\sum_{i=0}^{n-1} \left( \gamma (u_{i+1}-u_{i}) - W(i,i+1) e^{u_{i+1}-u_i} - W(i+1,i) e^{u_{i}-u_{i+1}}\right) \right) \prod_{i=1}^{n} \dd u_i.
$$
From this expression, it is clear that the sequence of random variables $V_i= (U_{i+1}-U_i)$ are independent with distribution on $\R$ proportional to:
$$
\exp\left(\gamma v - W(i,i+1) e^{v} - W(i+1,i) e^{-v} \right) \dd v.
$$
Since it only depends on~$i$ and not on the graph~$\cG_n$, it implies that we can define the infinite sequence of random variable $(V_i)_{i\in \NN}$ with that distribution and $(U_i)_{i\in \NN}$ by $U_0=0$ and $U_k=\sum_{i=0}^{k-1} V_i$. By construction, $(U_0, \ldots U_n)$ has distribution  $\nu_\gamma^{W, \cG_{n}}$.
	\end{proof}

Consider now $(U_i)_{i\in \NN}$ as defined above, and define $W^U$ as before by $W_{i,i+1}^U = W_{i,i+1} e^{U_{i+1}-U_i}$ and $W_{i+1,i}^U = W_{i+1,i} e^{U_{i}-U_{i+1}}$. By the previous lemma and Proposition~\ref{THmixing}, we have that $W^U$ has the same law as $W$. We define $\Gamma_n^U$ and $S_n^U$ from $W^U$ as $\Gamma_n$ and $S_n$ were defined  from $W$. Hence $((\Gamma_n^U)_{n\ge 1},(S_n^U)_{n\ge 1})$  has the same law as $((\Gamma_n)_{n\ge 1}, (S_n)_{n\ge 1})$.

We then define,
$$
\beta(0)=W_{0,1}e^{U_1-U_0}, \;\;\; \beta(i)=W_{i,i+1} e^{U_{i+1}-U_i}+ W_{i,i-1} e^{U_{i-1}-U_i}, \;\;\; \forall i\ge 1.
$$
If $\beta_n(k)$, $0\le k\le n$ is defined as in Corollary~\ref{THmixing_cor} for the graph $\cG_n$, with the above coupling of random variables $U$, we see that:
\begin{align}\label{beta_n=beta}
\beta_n(k)=\beta (k), \quad \forall 0 \pp k \pp n-1.
\end{align}
	
	\begin{lemma}
	\label{THMarkov}
	The distribution of $S^U_n$ conditionally to $(\Gamma^U_1, \ldots, \Gamma^U_n)$ is equal to the distribution of $S^U_n$ conditionally to $\Gamma^U_n$.
	\end{lemma}
	
	\begin{proof}
From the definition of $\omega^U$ and $\Gamma_i^U$, we have,	
for every $1\le i\le n$:
\begin{align*}
 \sigma(\Gamma^U_i) &\subseteq  \sigma(\ens{W(j,j\pm 1), 1 \pp j \pp i-1}, \beta_i(1), \ldots, \beta_i(i-1))\\
 &=\sigma(\ens{W(j,j\pm 1), 1 \pp j \pp i-1}, \beta(1), \ldots, \beta(i-1))
 \\
 &=\sigma(\ens{W(j,j\pm 1), 1 \pp j \pp i-1}, \beta_n(1), \ldots, \beta_n(i-1))
 \end{align*}
using the above coupling and \eqref{beta_n=beta}.
Therefore: for every $n \pg 1$,
	$$ \sigma(\Gamma^U_1, \ldots, \Gamma^U_n) \subseteq  \sigma(\ens{W(j,j\pm 1), 1 \pp j \pp n-1}, \beta_n(1), \ldots, \beta_n(n-1)). $$

	We have already noted in the proof of Lemma~\ref{THlaw} that the distribution of $S^U_n$ knowing $\Gamma^U_n$ is in fact the distribution of $S_n$ knowing $W$ and $(\beta_n(k))_{k \not\in \{0,n\}}$. Hence we have
\begin{align*}
\hbox{Law}\left( S^U_n \vert \Gamma^U_1, \ldots, \Gamma^U_n\right)
&=\hbox{Law}\left( S^U_n \vert \ens{W(j,j\pm 1), 1 \pp j \pp k-1}, \beta_n(1), \ldots, \beta_n(n-1)\right)\\
&=\hbox{Law}\left( S^U_n \vert \Gamma_n^U\right),
\end{align*}
which is what we need to prove.	\end{proof}
	
The conjunction of Lemma~\ref{THMarkov} and Lemma~\ref{THrecursion} proves that $(\Gamma^U_n)_{n \pg 1}$ is a Markov chain. As we noted that $(\Gamma_n^U)_{n\ge 1}$ and $(\Gamma_n)_{n\ge 1}$ have same distribution, we can conclude the proof.
%
%
	
\begin{remark} In comparison with Matsumoto-Yor property for Brownian motion, our discrete counterpart lacks a result about the strict inclusion of $(\cF_n)_{n \pg 1}$, the filtration adapted to $(\Gamma_n)_{n \pg 1}$, in $(\cS_n)_{n \pg 1}$, the filtration adapted to $(S_n)_{n \pg 1}$. However, our proof entails such a result when $W$ is fixed. In this case, $(\cF^U_n)_{n \pg 1}$, the filtration adapted to $(\Gamma^U_n)_{n \pg 1}$, is included in $\sigma(\beta(1), \ldots, \beta(n-1))$, as stated in the previous proof; on the other hand, $(\cS^U_n)_{n \pg 1}$, the filtration adapted to $(S^U_n)_{n \pg 1}$, is included in ${\sigma(U_1, \ldots, U_n)} = {\sigma(\beta(0), \ldots, \beta(n-1), \beta_n(n))}$. But note that by construction, $\beta_n(n)$ is redundant in this last sigma-field: actually, we have $(S^U_n)_{n \pg 1} \subseteq {\sigma(\beta(0), \ldots, \beta(n-1))}$.
Thus, when $W$ is deterministic, the difference between $(\cF^U_n)_{n \pg 1}$ and $(\cS^U_n)_{n \pg 1}$ lies in the knowledge of $\beta(0)$, the first of the $\beta$'s. This can be compared to the case of Brownian motion, where the missing information is the infinitesimal behaviour of the Brownian motion at its starting point.
\end{remark}	
	
\section{Proof of the non-existence of the invariant measure viewed from the particle in the recurrent case}
\label{recurrence}

This section contains the proof of Theorem~\ref{THrecurrence} and of the first part of Theorem~\ref{THaccelerated}.

	\subsection{Case of the original random walk}

We consider the random walk $X$ in Dirichlet environment $\PP$ over $\ZZ^2$, and assume that $X$ is almost surely recurrent. We aim at proving that $X$ cannot admit an absolutely continuous invariant measure from the point of view of the particle. We proceed by contradiction, and assume the existence of such an invariant measure $\QQ$, absolutely continuous with respect to $\PP$. As $\PP$ is invariant by translation over $\ZZ^2$, we deduce that $\pi^\omega=(\pi^\omega(x))_{x\in \ZZ^2}$ defined by
$$
\pi^\omega(x)={d\QQ\over d\PP}(\tau_x\omega),
$$
is an invariant measure for the quenched random walk $X$ for $\PP$-a.s. $\omega \in \Omega$. Indeed, if $f$ is positive measurable test function on $\Omega$, we have for all $x\in \ZZ^2$,
\begin{align*}
\int \sum_{i=1}^4 \pi^\w(x-e_i) \w(x-e_i,x) f(\tau_x \w) d\PP(\w)&= \int \sum_{i=1}^4 \pi^\w(0) \w(0,e_i) f(\tau_{e_i} \w) d\PP(\w)
\\
&=\int \cR f d\QQ(\w)
=\int f(\w) d\QQ(\w)
\\
&=\int \pi^\w(0) f(\w) d\PP(\w)= \int \pi^\w(x) f(\tau_x \w) d\PP(\w)
\end{align*}

We will successively show that, in that case, the ratio $\frac{\pi^\omega(j_0)}{\pi^\omega(i_0)}$ is small with high probability and big with high probability too, hence contradiction.

\begin{proof}[Proof of an upper bound]
On the one hand, note that, because of the invariance by translation of $\PP$, all the $\pi^\omega(x)$, $x \in \ZZ^2$ have the same distribution.
This entails that $\frac{\pi^\omega(j_0)}{\pi^\omega(i_0)}$, with high probability, cannot take too big values. Indeed, let $A, \eta > 0$; we have:
\begin{align*}
\PP\left(\frac{\pi^\omega(j_0)}{\pi^\omega(i_0)} \pg A \right) &= \PP\left(\pi^\omega(j_0) \pg A\pi^\omega(i_0), \pi^\omega(i_0) \pp \eta \right) + \PP\left(\pi^\omega(j_0) \pg A\pi^\omega(i_0), \pi^\omega(i_0) \pg \eta \right) \\
&\pp \PP\left(\pi^\omega(i_0) \pp \eta \right) + \PP\left(\pi^\omega(j_0) \pg A \eta \right) \\
&\pp \PP(\pi^\omega(0) \not\in [\eta, A\eta])
\end{align*}
The latter bound tends, when $A \to +\infty$, to $\PP(\pi^\omega(0) \le \eta)$. By ergodicity, $\PP(\pi^\omega(0) = 0) = 0$ (cf. \cite{Sznitman_ten_lectures}). Let us so fix $\eta > 0$ such that $\PP(\pi^\omega(0) < \eta) \pp \frac{1}{16}$ and $A > 0$ such that
$\PP(\pi^\omega(0) > A\eta)\le \frac{1}{16}$,
then: 
$$ \PP\left(\frac{\pi^\omega(j_0)}{\pi^\omega(i_0)} \pg A \right) \pp \frac{1}{8}. $$
\end{proof}

\begin{proof}[Proof of a lower bound]
On the other hand, since the quenched random walk $X$ is assumed to be recurrent, we know that the invariant measure is unique up to a multiplicative constant and that
$$
{\pi^\w(j_0)\over \pi^\w(i_0)}= \frac{\sP_{i_0}^{\omega, \ZZ^2} \left(H_{j_0} < H_{i_0}^+\right)}{\sP_{j_0}^{\omega, \ZZ^2} \left(H_{i_0} < H_{j_0}^+\right)}.
$$  
We can approximate the ratios $\frac{\pi^\omega(j_0)}{\pi^\omega(i_0)}$ of the invariant measure over $\ZZ^2$ by similar ratios of the invariant measure over a large torus. Indeed, let $\TT_N=(\ZZ/2N\ZZ)^2 - (N,N)$ the torus (with periodic boundary condition for the environment $\w$) and $\pi_N^\w$ be the associated invariant measure. We have for $i_0$ and $j_0$ in $\ZZ^2$:
$$ {\pi_N^\w(j_0)\over \pi_N^\w(i_0)}= \frac{\sP_{i_0}^{\omega, \TT_N} \left(H_{j_0} < H_{i_0}^+\right)}{\sP_{j_0}^{\omega, \TT_N} \left(H_{i_0} < H_{j_0}^+\right)}
\xrightarrow[N\to\infty]{\text{a.s.}}
{\pi^\w(j_0)\over \pi^\w(i_0)}=\frac{\sP_{i_0}^{\omega, \ZZ^2} \left(H_{j_0} < H_{i_0}^+\right)}{\sP_{j_0}^{\omega, \ZZ^2} \left(H_{i_0} < H_{j_0}^+\right)},$$
this limit being easily obtained by coupling the environments on the torus with the periodic environments on the plane. 
Therefore, for every $A > 0$: 
$$  \PP^{\ZZ^2}\left(\frac{\pi^\omega(j_0)}{\pi^\omega(i_0)} \pg A\right) = \lim_{N\to \infty} \PP^{\TT_N}\left(\frac{\pi^\omega_N(j_0)}{\pi^\omega_N(i_0)} \pg A\right). $$

Let us work out the latter probability, of which we know that it is linked, by Lemma~\ref{THidentity} applied to the torus, to the known variable~$S$. We have, for every $A, \eta > 0$:
\begin{align*}
\PP^{\TT_N}(e^{2S} \pg A)
&= \PP^{\TT_N}\left(\frac{\beta_{i_0}}{\beta_{j_0}} \frac{\pi^\omega_N(j_0)}{\pi^\omega_N(i_0)} \pg A\right) \\
&= \PP^{\TT_N}\left(\frac{\pi^\omega_N(j_0)}{\pi^\omega_N(i_0)} \pg A \frac{\beta_{j_0}}{\beta_{i_0}}, \frac{\beta_{j_0}}{\beta_{i_0}} \pg \eta\right) + \PP^{\TT_N}\left(\frac{\pi^\omega_N(j_0)}{\pi^\omega_N(i_0)} \pg A \frac{\beta_{j_0}}{\beta_{i_0}}, \frac{\beta_{j_0}}{\beta_{i_0}} \pp \eta\right) \\
&\pp \PP^{\TT_N}\left(\frac{\pi^\omega_N(j_0)}{\pi^\omega_N(i_0)} \pg A \eta\right) \PP^{\TT_N}\left(\frac{\beta_{j_0}}{\beta_{i_0}} \pg \eta\right) + \PP^{\TT_N}\left(\frac{\beta_{j_0}}{\beta_{i_0}} \pp \eta\right)
\end{align*}
where we used in the last line the independence of the $\beta$'s from $\omega$. We get:
$$ \PP^{\TT_N}\left(\frac{\pi^\omega_N(j_0)}{\pi^\omega_N(i_0)} \pg A\right) \pg \frac{\PP^{\TT_N}\left(e^{2S} \pg \frac{A}{\eta}\right) - \PP^{\TT_N}\left(\frac{\beta_{j_0}}{\beta_{i_0}} \pg \eta\right)}{1 - \PP^{\TT_N}\left(\frac{\beta_{j_0}}{\beta_{i_0}} \pg \eta\right)}.$$
Moreover, for every $\delta > 0$:
$$ \PP^{\TT_N}\left(e^{2S} \pg \frac{A}{\eta}\right) \pg \PP^{\TT_N}\left(\left. e^{2S} \pg \frac{A}{\eta} \right\vert \Gamma < \delta \right) \PP^{\TT_N}(\Gamma < \delta). $$
Finally, for every $A, \eta, \delta > 0$:
$$ \PP^{\TT_N}\left(\frac{\pi^\omega_N(j_0)}{\pi^\omega_N(i_0)} \pg A\right) \pg \frac{\PP^{\TT_N}\left(\left. e^{2S} \pg \frac{A}{\eta} \right\vert \Gamma < \delta \right) \PP^{\TT_N}(\Gamma < \delta) - \PP^{\TT_N}\left(\frac{\beta_{j_0}}{\beta_{i_0}} \pg \eta\right)}{1 - \PP^{\TT_N}\left(\frac{\beta_{j_0}}{\beta_{i_0}} \pg \eta\right)}. $$

By Laplace's method, we can show that, for every $\alpha > 0$, for $\delta$ small enough, we get:
\begin{equation}
\label{eqRec1}
\PP^{\TT_N}\left(\left. e^{2S} \pg \delta^{-\frac{1}{2}\alpha} \right\vert \Gamma < \delta\right) \pg \frac{1}{2}(1-\alpha).
\end{equation}

As we assume that the walk is recurrent, and thanks to Lemma~\ref{THidentity}, given some $\delta > 0$, we can choose $i_0$ and $j_0$ sufficiently far apart such that we have:
\begin{equation}
\label{eqRec2}
\PP^{\TT_N}(\Gamma < \delta) \pg 1 - \frac{1}{8}.
\end{equation}

Since the $\beta$ are continuous variables, there exists $\eta \pg 0$ such that we have:
\begin{equation}
\label{eqRec3}
\PP^{\TT_N}\left(\frac{\beta_{j_0}}{\beta_{i_0}} \pg \eta\right) \pg \frac{1}{8}.
\end{equation}

Fix the parameters as follows: first choose $\eta > 0$ such that \ref{eqRec3} is satisfied; then choose $\alpha < \frac{1}{8}$ and $\delta > 0$ small enough such that \ref{eqRec1} is satisfied; fix $i_0$ and $j_0$ such that \ref{eqRec2} is satisfied for this choice of $\delta$, and set $A = \eta \delta^{-\frac{1}{2} \alpha}$. Up to increasing $\eta$, we can have $A$ as large as we want. We finally get, that, for every $A > 0$ large enough, for every $N$ large enough so that the torus contains $i_0$ and $j_0$:
$$ \PP^{\TT_N} \left(\frac{\pi^\omega_N(j_0)}{\pi^\omega_N(i_0)} \pg A\right) \pg \frac{\frac{1}{2}(1-\frac{1}{8})(1-\frac{1}{8}) - \frac{1}{8}}{1 - \frac{1}{8}} \pg \frac{1}{4}. $$
Taking $N \to \infty$ gives, for every $A > 0$:
$$ \PP\left(\frac{\pi^\omega(j_0)}{\pi^\omega(i_0)} \pg A \right) \pg \frac{1}{4}. $$
\end{proof}

The contradiction arising from these two bounds proves Theorem~\ref{THrecurrence}.

	\subsection{Case of the accelerated random walk}
	
	We now examine an accelerated random walk $Y$ in Dirichlet environment $\PP$ over $\ZZ^2$ with acceleration factor $\gamma$ (recall that $\gamma$ has \emph{finite} support~$\Lambda$), and assume that $Y$ is almost surely recurrent. We still proceed by contradiction, and assume that there exists some absolutely continuous invariant measure~$\tilde\QQ$ from the point of view of the particle.
	
	As in the original case, we deduce the existence of an invariant measure $\bar\pi^\omega$ for the accelerated random walk for almost all $\omega \in \Omega$; the measure $\pi^\omega = \gamma(\tau_\cdot \omega) \bar\pi^\omega$ would then be invariant for the original random walk (but with possibly infinite mass). The ratios $\frac{\pi^\omega(j_0)}{\pi^\omega(i_0)} = \frac{\gamma(\tau_{j_0}\omega) \bar\pi^\omega(j_0)}{\gamma(\tau_{i_0}\omega) \bar\pi^\omega(i_0)}$ are then subject to the same contradiction as in the original case, provided that $i_0$ and $j_0$ are sufficiently far apart so that the supports $\tau_{i_0}\Lambda$ and $\tau_{j_0}\Lambda$ of $\gamma(\tau_{i_0}\omega)$ and $\gamma(\tau_{j_0}\omega)$ are disjoint---in which case those two variables are independent. The conclusion is then the same.
	
\section{Proof of the existence of the invariant probability viewed from the particle under condition~$\mathbf{(T)}$}
\label{transience}

This section is dedicated to the proof of Theorem~\ref{THtransience}, and, in the last subsection, of the second part of Theorem~\ref{THaccelerated}. We consider the random walk in Dirichlet environment $\cD(\alpha_1, \ldots, \alpha_4)$ over the plane $\ZZ^2$, and assume that it satisfies condition~$\mathbf{(T')}$.

The strategy is to approximate the plane by a sequence of growing toruses, for each one of which, as finite graphs, we can always find an invariant measure. Our main concern is to establish the convergence of (a subsequence of) the sequence of those invariant measures. A difficulty arising from the approximation by torus is that condition~$\mathbf{(T')}$ looses part of its strength: indeed, it is possible for the walk on the torus to reach its starting point by going all around the torus without any backtrack. To prevent this event, we introduce a killing property at any vertex of the torus by adding a ``cemetery'' vertex~$\partial$.

We consider the graph $\TT_N^*$ obtained from the torus centred at~$0$ $\TT_N = (\ZZ/2N\ZZ)^2 - (N,N)$ by adding a vertex $\partial$ and edges from every vertex of $\TT_N$ to $\partial$ and conversely. The graph $\TT_N^*$ and  is endowed with translation operators $(\tau_y)_{y\in\TT_N}$, by setting conventionally $y + \partial = \partial$ for every $y \in \TT_N$. An infinite graph $\TT^*$ is obtained from $\ZZ^2$ with the same procedure.

We denote by $\Omega_N$, $\Omega_N^*$ and $\Omega^*$ the space of environments on $\TT_N$, $\TT_N^*$ and $\TT^*$ respectively. We denote $\PP^{(\alpha)}_N$ (respectively $\PP^{(\alpha)}$) the Dirichlet environment on $\TT_N$ (respectively on $\ZZ$). Given a small $\epsilon > 0$, we deduce from them a Dirichlet environment $\PP^{(\alpha,\epsilon)}_N$ on $\TT_N^*$ (respectively on $\PP^{(\alpha,\epsilon)}$ on $\TT^*$) by keeping $\alpha$ on $\TT_N$ and endowing with a weight $\epsilon$ the edges going to or coming from $\partial$. We can extend by periodicity a Dirichlet environment $\PP^{(\alpha,\epsilon)}_N$ on $\TT_N^*$ in a Dirichlet environment on $\TT^*$.

	\subsection{From the plane to the torus}
	
We first use condition~$\mathbf{(T')}$ to estimate the probability for the walk starting from~0 to hit~0 again before being killed or before crossing the whole torus. The proof is divided in two steps: first, we use condition~$\mathbf{(T')}$ to control the probability of returning to~$0$ (or, equivalently, the average number of returns to~$0$) on~$\ZZ^2$; then, we transfer this estimation valid for $\ZZ^2$ to the torus $\TT_N^*$.

Let $N_0(T)$ be the number of visits to~$0$ before a time~$T$, and $N_0$ the total number of visits to~$0$.

\begin{lemma}
\label{THplane}
On $\ZZ^2$, under condition~$\mathbf{(T')}$, for all $\beta \in (1, \kappa)$, we have:
$$ \sE_0\left[{N_0}^\beta\right] < \infty. $$
\end{lemma}

\begin{proof}
Note that $N_0 = N_0(T_1)$, the number of visits to~$0$ before the first renewal time~$T_1$. By condition~$\mathbf{(T')}$, we can assume that this renewal time occurs before the walk goes too far away from~$0$: for $x > 0$, let $\ell(x) = \left(\frac{\kappa}{c} \ln(x)\right)^2$ and $U(x)$ be a box centred at~$0$ with width $2\ell(x)$; by condition~$\mathbf{(T')}$, there exists a finite constant~$C$ such that:
$$ \sP_0(T_1 > \bar H_{U(x)}) \pp Ce^{-c\ell(x)^{1/2}} \pp C x^{-\kappa}. $$

On the event that the first renewal time occurs before exiting $U(x)$, we can use a very similar proof to lemma~16 in \cite{perrel2024limit} to show, for every $1 < \beta < \kappa$, the existence of two constants $c_\beta$ and $\gamma_\beta$ such that $ \sE_0\left[(N_0)^\beta \Ind{T_1 < \bar H_{U(x)}}\right] \pp c_\beta \ell(x)^{\gamma_\beta}$. We do not detail this proof here, because we will detail it in the accelerated case: cf. Lemma~\ref{THplaneBIS}. This ensures, for every $0 < \beta < \kappa$, the existence of two constants $c_\beta$ and $\gamma_\beta$ such that:
$$ \sP_0(N_0 \pg x, T_1 < \bar H_{U(x)}) \pp c_\beta x^{-\beta} \ln(x)^{\gamma_\beta}.$$

Combining those two inequalities, we finally get that: for every $\beta \in (1, \kappa)$,
$$ \sP_0(N_0 \pg x) = o\left(x^{-\beta}\right). $$
This is sufficient to complete the proof.
\end{proof}

\begin{lemma}
\label{THvisits}
On $\TT_N^*$, for all $\beta \in (1, \kappa)$ there exists a constant $C_\beta \in (0, +\infty)$ such that, for $N$ large enough:
$$ \EE^{(\alpha, \epsilon)}_N\left[\sP_0^\omega(H_\partial < H_0^+)^{-\beta}\right] \pp \sE_0\left[{N_0(H_\partial)}^\beta\right] < C_\beta. $$
\end{lemma}

\begin{proof}
On $\TT_N^*$, for every $\omega \in \Omega_N^*$, $N_0(H_\partial)$ follows a geometric distribution with parameter $\sP_0^\omega(H_\partial < H_0^+)$. Therefore, for any $\beta > 1$, by Jensen's inequality:
\begin{align*}
\EE_N\left[\sP_0^\omega(H_\partial < H_0^+)^{-\beta}\right]
= \EE_N\left[\sE_0^\omega[N_0(H_\partial)]^{\beta}\right]
\pp \sE_0\left[{N_0(H_\partial)}^\beta\right].
\end{align*}

The random variable $N_0(H_\partial)$, which counts the number of returns to~$0$ before hitting~$\partial$, can be stochastically dominated under $\sP_0^\omega$ by a random sum of i.i.d. random variables $\sum_{i=1}^Y N'_i$ as follows. Set $\Theta_0=0$; denote by $\Theta_1$ the first hitting time of the boundary of~$\TT_N$ (\textit{i.~e.} the set of vertices at distance~$N$ of the origin union $\{\partial\}$) and $N'_1$ the number of visits to~$0$ before~$\Theta_1$; recursively, denote by $\Theta_{i+1}$ the first time the walk reaches the boundary of~$\TT_N$ after the first visit to~$0$ following $\Theta_i$, and $N'_{i+1}$ the number of visits to~$0$ between $\Theta_i$ and~$\Theta_{i+1}$. Obviously, the $N'_i$'s are stochastically dominated by $N_0$, the number of visits to~$0$ in~$\ZZ^2$: so the $N'_i$'s are independent copies of $N_0$ and $Y$ is the number of $\Theta_i$'s occurred before $H_\partial$ included. As a consequence, by a classical convexity inequality, for any $\beta > 1$:
\begin{align*}
\sE_0\left[{N_0(H_\partial)}^\beta\right]
\pp \sE_0\left[\left(\sum_{i=1}^Y N'_i\right)^\beta\right]
&\pp \sE_0\left[\sum_{k\pg 0} \left(\sum_{i=1}^k N'_i\right)^\beta \Ind{Y=k}\right] \\
&\pp \sum_{k\pg 0} k^{\beta-1} \sum_{i=1}^k \sE_0\left[{N'_i}^\beta \Ind{Y=k}\right]
\end{align*}
and, by H\"older's inequality, for every $\frac{1}{p} + \frac{1}{q} = 1$:
\begin{align*}
\sE_0\left[{N_0(H_\partial)}^\beta\right]
\pp \sum_{k\pg 0} k^{\beta-1} \sum_{i=1}^k \sE_0\left[{N'_i}^{p\beta}\right]^{\frac{1}{p}} \sP_0(Y=k)^\frac{1}{q}.
\end{align*}

Fix $\beta \in (1, \kappa)$ and $p > 1$ such that $p\beta < \kappa$. Thanks to the previous lemma, we already know that there exists a finite constant $C_{p\beta} > 0$ such that, for every $i \pg 1$, $\sE_0\left[{N'_i}^{p\beta}\right] \pp C_{p\beta}$. Hence:
\begin{align*}
\sE_0\left[{N_0(H_\partial)}^\beta\right]
\pp {C_{p\beta}}^\frac{1}{p} \sum_{k\pg 0} k^{\beta} \sP_0(Y=k)^\frac{1}{q}.
\end{align*}

There remains to determine an equivalent of $\sP_0(Y=k)$ when $k \to \infty$. Note that, if, starting from $0$, the walk hits the boundary of $\TT_N$ somewhere else than at $\partial$ (as is the case at $\Theta_1, \ldots, \Theta_{Y-1}$), this means that the walks meets a set $\cI$ of at least $N$ vertices of $\TT_N$ (for example $0$, the first vertex encountered at distance~$1$ of the origin, the first vertex at distance~$2$, and so on and so forth until the first vertex at distance $N$), and that, at each one of those vertices, the walks chooses to go to a neighbour in~$\TT_N$ rather than to~$\partial$---the probability of this event following a beta distribution $\beta\left(\sum_{i=1}^4 \alpha_i, \epsilon\right)$ under~$\sP_0$. Therefore, if we arbitrarily enumerate the vertices of $\TT_N$ and denote by $\Pi_j$ the probability for the walk at vertex~$j$ to go to a neighbour of~$j$ in~$\TT_N$ rather than to~$0$, for every $k \pg 1$, we can partition the event $\ens{Y = k}$ according to the values $I_1, \ldots, I_{k-1}$ of the sets $\cI_1, \ldots, \cI_{k-1}$ of these $N$~vertices encountered by the walk during $[\Theta_0,\Theta_1], \ldots, [\Theta_{k-2},\Theta_{k-1}]$:
$$ \sP_0(Y = k) \pp \sum_{\substack{I_1, \ldots, I_k \subset \NN \tq \\ \#I_1 = \ldots = \#I_k = N}} \sE_0\left[\prod_{i=1}^{k-1} \prod_{j \in I_i} \Pi_i \right] \sP_0\left(\cI_1 = I_1, \ldots, \cI_k = I_k\right). $$
By Lemma~\ref{THappendix} (in appendix), we have, whatever the values of $I_1, \ldots, I_k$:
$$\forall 1 \pp i \pp k, \quad \sE_0\left[\prod_{i=1}^{k-1} \prod_{j \in I_i} \Pi_i \right] \pp \sE_0\left[{\Pi_1}^{k-1}\right]^{N},$$
so that:
\begin{align*}
\sP_0(Y = k) \pp \sE_0\left[{\Pi_1}^{k-1}\right]^{N}
&= \left(\frac{\Gamma\left(\sum_{i=1}^{4}\alpha_i + \epsilon\right)\Gamma\left(\sum_{i=1}^{4}\alpha_i + k\right)}{\Gamma\left(\sum_{i=1}^{4}\alpha_i\right)\Gamma\left(\sum_{i=1}^{4}\alpha_i + \epsilon+k\right)}\right)^{N} \\
&= \left(\prod_{i=0}^{k-1}\frac{\sum_{i=1}^{4}\alpha_i+i}{\sum_{i=1}^{4}\alpha_i+\epsilon+i}\right)^{N} \\
&\underset{k\to\infty}{\sim} k^{-N\epsilon}.
\end{align*}

Finally, remind that:
\begin{align*}
\sE_0\left[{N_0(H_\partial)}^\beta\right]
\pp {C_{p\beta}}^\frac{1}{p} \sum_{k\pg 0} k^{\beta} \sP_0(Y=k)^\frac{1}{q}.
\end{align*}
The summand is equivalent to $k^{\beta-\frac{N\epsilon}{q}}$ and this sum is finite as soon as $N > \frac{q\beta}{\epsilon}$.
\end{proof}

	\subsection{Uniform bounds for invariant measures on growing torus}

Since $\TT_N^*$ is finite and strongly connected, and since $\PP^{(\alpha)}\text{-a.s.}$ the $W$'s are positive, there exists a unique invariant measure $\pi^\omega_N$ for the environment~$\omega$.

We define:
$$ \dd \QQ^{(\alpha, \epsilon)}_N = \frac{\pi^\omega_N(0)/\tilde\pi^\omega_N(\partial)}{\EE_N\left[\pi^\omega_N(0)/\tilde\pi^\omega_N(\partial)\right]} \dd \PP_N^{(\alpha, \epsilon)} $$
where $\tilde\pi^\omega_N(\partial) = \frac{\pi^\omega_N(\partial)}{\beta_\partial}$ and $\beta_\partial$ is a random variable following $\Gamma(4N^2\epsilon,1)$ independent from~$\omega$.

Because the expectation at the denominator in the definition of $\QQ_N^{(\alpha, \epsilon)}$ is possibly infinite, we must first check that $\QQ_N^{(\alpha, \epsilon)}$ is non-null. We can even show more precisely that the sequences of functions $\left(\frac{\dd\QQ_N^{(\alpha, \epsilon)}}{\dd\PP_N^{(\alpha, \epsilon)}}\right)_{N\pg1}$ is uniformly lower- and upper-bounded in some $L^a, {a > 1}$:

\begin{lemma}
\label{THbounded}
Assume that $\kappa > 1$. For all $a \in (1, \kappa)$ there exist finite constants $c, C > 0$ such that, for $N$ large enough:
$$ c \pp \EE_N^{(\alpha, \epsilon)}\left[ \left( \frac{\pi^\omega_N(0)}{\tilde\pi_\omega^N(\partial)} \right)^a \right] \pp C. $$
and
$$ c \pp \EE_N^{(\alpha, \epsilon)}\left[ \left\| \frac{\dd\QQ_N^{(\alpha, \epsilon)}}{\dd\PP_N^{(\alpha, \epsilon)}} \right\|^a \right] \pp C. $$
\end{lemma}

\begin{proof}
The second inequality is an immediate consequence of the first one. Let $a \in (1, \kappa)$.

\paragraph{\textbf{Upper-bound.}} We have, by H\"older's inequality, for every $\frac{1}{p} + \frac{1}{q} = 1$:
\begin{align*}
\EE_N^{(\alpha, \epsilon)}\left[ \left( \frac{\pi^\omega_N(0)}{\tilde\pi_\omega^N(\partial)} \right)^a \right]
&= \EE_N^{(\alpha, \epsilon)}\left[ \left( \frac{\beta_\partial \sP_\partial^\omega(H_0 < H_\partial^+)}{\sP_0(H_\partial < H_0^+)} \right)^a \right] \\
&\pp \EE_N^{(\alpha, \epsilon)}\left[\sP_0(H_\partial < H_0^+)^{-ap}\right]^{\frac{1}{p}} \EE_N^{(\alpha, \epsilon)}\left[\left(\beta_\partial \sP_\partial^\omega(H_0 < H_\partial^+)\right)^{aq}\right]^{\frac{1}{q}}.
\end{align*}
Note that, by lemma~\ref{THidentity}, $\beta_\partial \sP_\partial^\omega(H_0 < H_\partial^+)$ and $\beta_0 \sP_0^\omega(H_\partial < H_0^+)$, where $\beta_0$ is a random variable following $\Gamma(\alpha(0),1)$ independently from~$\omega$, are identically distributed. Hence:
\begin{align*}
\EE_N^{(\alpha, \epsilon)}\left[ \left\| \frac{\pi^\omega_N(0)}{\tilde\pi_\omega^N(\partial)} \right\|^a \right]
&\pp \EE_N^{(\alpha, \epsilon)}\left[\sP_0(H_\partial < H_0^+)^{-ap}\right]^{\frac{1}{p}} \EE_N^{(\alpha, \epsilon)}\left[\left(\beta_0 \sP_0^\omega(H_\partial < H_0^+)\right)^{aq}\right]^{\frac{1}{q}} \\
&\pp \EE_N^{(\alpha, \epsilon)}\left[\sP_0(H_\partial < H_0^+)^{-ap}\right]^{\frac{1}{p}} \EE_N^{(\alpha, \epsilon)}\left[{\beta_0}^{aq}\right]^{\frac{1}{q}}.
\end{align*}
By lemma~\ref{THvisits}, for all $a \in (1, \kappa)$ and $p > 1$ such that $ap<\kappa$ the first factor is bounded by a constant independent of $N$; the second factor is upper bounded for every $a > 1$ and $q > 1$.
 
\paragraph{\textbf{Lower-bound.}} Fix $a\in (1,\kappa)$. We have, by Jensen's inequality, for every $r > 0$:
\begin{align*}
\EE_N^{(\alpha, \epsilon)}\left[ \left\| \frac{\pi^\omega_N(0)}{\tilde\pi_\omega^N(\partial)} \right\|^a \right]
= \EE_N^{(\alpha, \epsilon)}\left[ \left( \beta_0 \frac{\tilde\pi^\omega_N(0)}{\tilde\pi_\omega^N(\partial)} \right)^a \right]
&= \EE_N^{(\alpha, \epsilon)}\left[{\beta_0}^a\right] \EE_N^{(\alpha, \epsilon)}\left[ \left( \frac{\tilde\pi^\omega_N(0)}{\tilde\pi_\omega^N(\partial)} \right)^a \right] \\
&\pg \EE_N^{(\alpha, \epsilon)}\left[{\beta_0}^a\right] \frac{1}{\EE_N^{(\alpha, \epsilon)}\left[ \left( \frac{\tilde\pi^\omega_N(\partial)}{\tilde\pi_\omega^N(\partial)} \right)^{ar} \right]^\frac{1}{r}}.
\end{align*}
In an analogous way as in the upper-bound, for every $\frac{1}{p} + \frac{1}{q} = 1$, we have:
\begin{align*}
\EE_N^{(\alpha, \epsilon)}\left[ \left( \frac{\pi^\omega_N(\partial)}{\tilde\pi_\omega^N(0)} \right)^{ar} \right]
&= \EE_N^{(\alpha, \epsilon)}\left[ \left( \frac{\beta_0 \sP_0^\omega(H_\partial < H_0^+)}{\beta_\partial \sP_\partial(H_0 < H_\partial^+)} \right)^a \right] \\
&\pp \EE_N^{(\alpha, \epsilon)}\left[\left(\beta_\partial \sP_\partial(H_0 < H_\partial^+)\right)^{-arp}\right]^{\frac{1}{p}} \EE_N^{(\alpha, \epsilon)}\left[\left(\beta_0 \sP_0^\omega(H_\partial < H_0^+)\right)^{arq}\right]^{\frac{1}{q}} \\
&\pp \EE_N^{(\alpha, \epsilon)}\left[{\beta_0}^{-arp}\right] \EE_N\left[\sP_0(H_\partial < H_0^+)^{-arp}\right]^{\frac{1}{p}} \EE_N^{(\alpha, \epsilon)}\left[{\beta_0}^{arq}\right]^{\frac{1}{q}}
\end{align*}
The first factor is finite as soon as $arp < \alpha(0)$; as in the upper-bound, there exists some 
$p > 1$ and $r > 0$ such that the second factor is upper-bounded; the third factor is always finite. We can so choose $r > 0$ small enough such that $\EE_N\left[ \left( \frac{\pi^\omega_N(\partial)}{\tilde\pi_\omega^N(0)} \right)^{ar} \right]$ is finite. This choice entails the desired lower-bound.
\end{proof}

The measure $\QQ_N^{(\alpha, \epsilon)}$ being non-null, obviously:

\begin{lemma}
\label{THprobameasure}
For every $N \pg 1$, $\QQ_N^{(\alpha, \epsilon)}$ is a probability measure over~$\Omega_N^*$.
\end{lemma}

We know examine the question of the invariance of $\QQ_N^{(\alpha, \epsilon)}$ viewed from the particle.

\begin{lemma}
\label{THapproximatedInvariance}
For every $N \pg 1$, $\QQ_N^{(\alpha, \epsilon)}$ satisfies the following ``approximated invariance'' equation: for every $f:\Omega \to \RR_+ \text{ measurable}$,
\begin{align*}
\int_\Omega f(\omega) \dd\QQ_N^{(\alpha, \epsilon)}
= \sum_{x \neq \partial \tq x \rightarrow 0} \int_\Omega f(\tau_{x}\omega) \omega(0,x) \dd\QQ_N^{(\alpha, \epsilon)}
+ \int_\Omega f(\omega) \omega(\partial, 0) \frac{\beta_\partial}{\EE_N^{(\alpha, \epsilon)}\left[\pi^\omega_N(0)/\tilde\pi^\omega_N(\partial)\right]} \dd\PP_N^{(\alpha, \epsilon)}.
\end{align*}
\end{lemma}

\begin{proof}
Since $\pi_N^\omega$ is invariant for~$\omega$, we have:
$$\pi^\omega_N(0) = \sum_{x \in \TT_N^* \tq x \rightarrow 0} \omega(x,0) \pi^\omega_N(x).$$
Let $f:\Omega \to \RR_+$ be a non-negative measurable function. We integrate the previous equality:
$$\int_\Omega f(\omega) \dd\QQ_N^{(\alpha, \epsilon)} = \sum_{x \tq x \rightarrow 0} \int_\Omega f(\omega) \omega(x,0) \frac{\pi^\omega_N(x)/\tilde\pi^\omega_N(\partial)}{\EE_N^{(\alpha, \epsilon)}\left[\pi^\omega_N(0)/\tilde\pi^\omega_N(\partial)\right]} \dd\PP_N^{(\alpha, \epsilon)}.$$
By invariance by translation of $\PP_N^{(\alpha, \epsilon)}$, and because, by uniqueness of the invariant measure for $\omega$, we have $\forall (u, v) \in \TT_N^*\times\TT_N, \pi^{\tau_v\omega}_N(u) = \pi^\omega_N(u+v)$:
\begin{align*}
\int_\Omega f(\omega) \dd\QQ_N^{(\alpha, \epsilon)} &= \sum_{x \neq \partial \tq x \rightarrow 0} \int_\Omega f(\tau_{-x}\omega) \omega(0,-x) \frac{\pi^{\tau_{-x}\omega}_N(x)/\tilde\pi^{\tau_{-x}\omega}_N(\partial)}{\EE_N\left[\pi^\omega_N(0)/\tilde\pi^\omega_N(\partial)\right]} \dd\PP_N^{(\alpha, \epsilon)} \\*
&\qquad + \int_\Omega f(\omega) \omega(\partial, 0) \frac{\pi^\omega_N(\partial)/\tilde\pi^\omega_N(\partial)}{\EE_N^{(\alpha, \epsilon)}\left[\pi^\omega_N(0)/\tilde\pi^\omega_N(\partial)\right]} \dd\PP_N^{(\alpha, \epsilon)} \\
&= \sum_{x \neq \partial \tq x \rightarrow 0} \int_\Omega f(\tau_{-x}\omega) \omega(0,-x) \frac{\pi^\omega_N(0)/\tilde\pi^\omega_N(\partial)}{\EE_N^{(\alpha, \epsilon)}\left[\pi^\omega_N(0)/\tilde\pi^\omega_N(\partial)\right]} \dd\PP_N^{(\alpha, \epsilon)} \\*
&\qquad + \int_\Omega f(\omega) \omega(\partial, 0) \frac{\beta_\partial}{\EE_N^{(\alpha, \epsilon)}\left[\pi^\omega_N(0)/\tilde\pi^\omega_N(\partial)\right]} \dd\PP_N^{(\alpha, \epsilon)}  \\
&= \sum_{x \neq \partial \tq x \rightarrow 0} \int_\Omega f(\tau_{x}\omega) \omega(0,x) \frac{\pi^{\omega}_N(0)/\tilde\pi^\omega_N(\partial)}{\EE_N^{(\alpha, \epsilon)}\left[\pi^\omega_N(0)/\tilde\pi^\omega_N(\partial)\right]} \dd\PP_N^{(\alpha, \epsilon)} \\*
&\qquad + \int_\Omega f(\omega) \omega(\partial, 0) \frac{\beta_\partial}{\EE_N^{(\alpha, \epsilon)}\left[\pi^\omega_N(0)/\tilde\pi^\omega_N(\partial)\right]} \dd\PP_N^{(\alpha, \epsilon)}\\
&= \sum_{x \neq \partial \tq x \rightarrow 0} \int_\Omega f(\tau_{x}\omega) \omega(0,x) \dd\QQ_N^{(\alpha, \epsilon)} \\*
&\qquad + \int_\Omega f(\omega) \omega(\partial, 0) \frac{\beta_\partial}{\EE_N^{(\alpha, \epsilon)}\left[\pi^\omega_N(0)/\tilde\pi^\omega_N(\partial)\right]} \dd\PP_N^{(\alpha, \epsilon)}.
\end{align*}
\end{proof}

Therefore, if we had $\beta_\partial \omega(\partial, 0)=0$ (\emph{i.e.} the second term of the right-hand side vanished), $\QQ_N^{(\alpha,\epsilon)}$ would be invariant from the point of view of the particle, because the last inequality is true for every non-negative measurable function~$f$. However, the second term vanishes only at the limit $\epsilon \to 0$:

\begin{lemma}
\label{THepsilonConvergence}
There exists a constant $C \in (0,\infty)$ such that, for every $N \pg 1$ and $f: \Omega \longrightarrow \RR$ bounded continuous function:
$$ \left\vert \int_\Omega f(\omega) \omega(\partial, 0) \frac{\beta_\partial}{\EE_N^{(\alpha, \epsilon)}\left[\pi^\omega_N(0)/\tilde\pi^\omega_N(\partial)\right]} \dd\PP_N^{(\alpha, \epsilon)}  \right\vert \pp C \|f\|_\infty \epsilon. $$
\end{lemma}

\begin{proof}
Remind that $\beta_\partial \omega(\partial, 0) = W(\partial, 0)$, follows the distribution $\Gamma(\epsilon, 1)$ under $\PP_N^{(\alpha, \epsilon)}$. Therefore,
$$ \EE_N^{(\alpha, \epsilon)}\left[ \beta_\partial \omega(\partial, 0)\right] = \epsilon. $$
From that, it readily follows that:
$$ \left\vert \int_\Omega f(\omega) \omega(\partial, 0) \frac{\beta_\partial}{\EE_N^{(\alpha, \epsilon)}\left[\pi^\omega_N(0)/\tilde\pi^\omega_N(\partial)\right]} \dd\PP_N^{(\alpha, \epsilon)} \right\vert \pp \frac{\|f\|_\infty}{\EE_N^{(\alpha, \epsilon)}\left[\pi^\omega_N(0)/\tilde\pi^\omega_N(\partial)\right]} \epsilon. $$
But we proved in Lemma~\ref{THbounded} that $\EE_N^{(\alpha, \epsilon)}\left[\pi^\omega_N(0)/\tilde\pi^\omega_N(\partial)\right]$ is uniformly bounded on $N$. This concludes the proof.
\end{proof}

	\subsection{Back from the torus to the plane}
	
	With the previous results of this section, the proof of Theorem~\ref{THtransience} is a double application of a classical argument found \textit{e.g.} in \cite{Sznitman_ten_lectures}.
	
	Fix $\epsilon > 0$. Extend $\PP_N^{(\alpha, \epsilon)}$ and $\QQ_N^{(\alpha, \epsilon)}$ as probability measures on $2N$-periodic environments of $\Omega^*$. When $N$ goes to $\infty$, it is obvious that $\PP_N^{(\alpha, \epsilon)}$ converges weakly to the probability measure $\PP^{(\alpha, \epsilon)}$ on $\Omega^*$.
	On the other hand, $\left(\QQ_N^{(\alpha,\epsilon)}\right)_{N\pg 1}$ is a sequence of probability measures. Since $\Omega^*$ is compact set (as a product of compact sets), so is the set of probability measures on $\Omega^*$.
	By compacity, up to extracting a subsequence, when $N$ goes to $\infty$, $\QQ_N^{(\alpha,\epsilon)}$ converges weakly to a probability measure $\QQ^{(\alpha,\epsilon)}$ on $\Omega^*$.
	Moreover, as $\left(\frac{\dd\PP_N^{(\alpha,\epsilon)}}{\dd\QQ_N^{(\alpha,\epsilon)}}\right)_{N \pg 1}$ is bounded (and note that it is uniformly bounded in $\epsilon$), $\QQ^{(\alpha,\epsilon)}$ is absolutely continuous with respect to $\PP^{(\alpha,\epsilon)}$.
	
	Now let $\epsilon$ converge to~$0$. Another application of this extraction argument for a decreasing sequence of $\epsilon$ to~$0$ gives the convergence of (a subsequence of) $\QQ^{(\alpha,\epsilon)}$ when $\epsilon$ goes to~$0$ to a limit $\QQ^{(\alpha)}$ over $\Omega$ absolutely continuous with respect to $\PP^{(\alpha)}$.

There remains to prove that $\QQ^{(\alpha)}$ is invariant from the point of view of the particle. This a direct consequence of Lemmas~\ref{THapproximatedInvariance} and~\ref{THepsilonConvergence}: given a bounded function~$f$, taking $N$ to $\infty$ in the approximated invariance equation of Lemma~\ref{THapproximatedInvariance} shows that $\QQ_{(\alpha,\epsilon)}$ satisfies the same equation; finally, thanks to Lemma~\ref{THepsilonConvergence}, $\epsilon \to 0$ makes the defect of invariance disappear: so $\QQ^{(\alpha)}$ is invariant.

Finally, the uniqueness of the invariant measure when there exists one is a classical result (see \cite{Sznitman_ten_lectures}).
	
	\subsection{Case of the accelerated random walk}
	\label{accelerationBis}
	
	As for the accelerated random walk under condition~$\mathbf{(T')}$, the existence of an acceleration factor for which there exists an absolutely continuous invariant probability from the point of view of the particle in the case $d=2$, even when $\kappa < 1$, comes from a similar reasoning for the case $d \pg 3$ of Bouchet in \cite{bouchet2013subballistic}.
	
	Because of non-uniform ellipticity of the environment, traps of finite size are created, \emph{i.e.} finite subsets in which the random walk spends an atypically large time. The strength of a possible trap ${S \subset V}$ can be defined as the order of the tail of the quenched Green function $G_S^\omega$ of the walk killed when exiting~$S$. Given a subset of vertices~$\Lambda \subset V$, we can measure how much $\Lambda$ delays the walk thanks to the strength~$\kappa(\Lambda)$ of the strongest finite trap containing~0 and not contained in~$\Lambda$. In a finite graph, Tournier proved in \cite{tournier2009integrability} that:
\begin{equation}
\kappa(\Lambda) = \min \ens{ \sum _{e \in \partial_+(S)}  \alpha _e; S \text{ connected set of vertices},\, \{0\} \subsetneq S \text{ and } \partial \Lambda \cap S \neq \emptyset }
\end{equation}
where we denote, for $S \subset V$:
$$\partial S = \ens{x \in S \tq \exists y \not\in S, x \sim y}, \qquad \partial_+ S = \ens{e \in E \tq \underline{e} \in S, \overline{e} \not\in S}. $$

	Note that, in $\ZZ^d$, the parameter $\kappa$ introduced in Theorem~\ref{THasymptotic1} turns out to be:
\begin{align*}
\kappa = \min \ens{\kappa(\Lambda), \Lambda \text{ finite connected set of vertices of } \ZZ^d} 
\end{align*}
In other words, in $\ZZ^d$, $\kappa$ is the minimal value of the $\kappa$-parameter of finite subsets, and undirected edges are the strongest finite traps of $\ZZ^d$. In fact, \cite{tournier2009integrability} has also proved that they are the strongest traps of $\ZZ^d$, stronger than any infinite subset of vertices. (See also \cite{Slonim2024} where the role of finite traps appear in a more subtle way).

	
	We use the local accelerating function Bouchet introduced in \cite{bouchet2013subballistic}:
$$ \gamma(\omega) = \frac{1}{\sum_{\sigma \in \Pi_\Lambda} \omega_\sigma}, $$
where, given a box~$\Lambda = [-n, n]^d$, the sum is extended to the set $\Pi_\Lambda$ of all directed paths $\sigma$ starting from $0$, ending outside $\Lambda$ and confined, except for the last vertex, in $\Lambda$. This acceleration counter the effect of all traps whose size is smaller than the size of~$\Lambda$, and where $\omega_\sigma = \prod_{i=1}^{|\sigma|} \omega(\sigma(i-1), \sigma(i))$ for every path~$\sigma$.

Our reasoning for the accelerated walk~$Y$ is very similar to the one for the original walk~$X$. On $\TT_N^*$, there exists a unique probability measure $\pi_N^\omega$ for the original random walk in environment $\omega$. Then the accelerated random walk admits an invariant measure $\bar\pi_N^\omega$ given by:
$$ \forall x \in \TT_N^*, \pi_N^\omega(x) = \gamma(\tau_x \omega) \bar\pi_N^\omega. $$

We introduce:
$$ \dd \bar\QQ_N^{(\alpha,\epsilon)} = \frac{\beta_\partial \bar\pi^\omega_N(0)/\pi^\omega_N(\partial)}{\EE_N^{(\alpha,\epsilon)}\left[\beta_\partial \bar\pi^\omega_N(0)/\pi^\omega_N(\partial)\right]} \dd \PP_N^{(\alpha,\epsilon)}. $$

As in the original case, $\bar\QQ_N$ is a sequence of probability measures which satisfies an approximated invariance equation with a extra term, similar to the one of Lemma~\ref{THapproximatedInvariance}. To control the extra term, we will prove that for all $a \in (1,\kappa(\Lambda))$, we can lower-bound and upper-bound the following quantity: 
$$ \EE_N^{(\alpha,\epsilon)} \left[\left(\beta_\partial \frac{\bar\pi_N^\omega(0)}{\pi_N^\omega(\partial)}\right)^a\right]. $$
The proof of Lemma~\ref{THbounded} is still valid in the accelerated setting, with the difference that we must adapt Lemma~\ref{THvisits} so as to control, instead of $\EE^{(\alpha,\epsilon)}_N\left[\sP_0(H_\partial < H_0^+)^{-ap}\right]$, the following quantity:
\begin{align*}
\EE^{(\alpha,\epsilon)}_N\left[\gamma(\omega)^{-ap} \sP_0(H_\partial < H_0^+)^{-ap}\right]
&\pp {\#\Pi_\Lambda}^{ap-1} \sum_{\sigma\in\Pi_\Lambda} \EE^{(\alpha,\epsilon)}_N\left[(\omega_\sigma)^{ap} \sP_0(H_\partial < H_0^+)^{-ap}\right] \\
&\pp {\#\Pi_\Lambda}^{ap-1} \sum_{\sigma\in\Pi_\Lambda} \sE_{N,0}^{(\alpha,\epsilon)}\left[(\omega_\sigma)^{ap} N_0(H_\partial)^{ap}\right].
\end{align*}
We thus want to bound ${\sE_{N,0}^{(\alpha,\epsilon)}}\left[(\omega_\sigma)^{ap} N_0(H_\partial)^{ap}\right]$ uniformly in~$\sigma$. Fix $\sigma \in \Pi_\Lambda$. Again, the proof of Lemma~\ref{THvisits} remains valid replacing $N_0$ by $\omega_\sigma N_0$. This, in turn, compells us to adapt Lemma~\ref{THplane}.

\begin{lemmabis}{THplane}
\label{THplaneBIS}
In $\ZZ^2$, fix $\Lambda = [-n, n]^2$. Under condition~$\mathbf{(T')}$, for every $\beta < \kappa(\Lambda)$, there exists a constant $C_\beta < \infty$ such that, for every $\sigma\in\Pi_\Lambda$:
$$ \sE_0^{(\alpha)}\left[{(\omega_\sigma N_0)}^\beta\right] < C_\beta. $$
\end{lemmabis}

\begin{proof}
Fix $\Lambda = [-n, n]^2$, $\beta < \kappa(\Lambda)$ and $\sigma \in \Pi_\Lambda$. As in Lemma~\ref{THplane}, by condition~$\mathbf{(T')}$, it suffices to show, for every $\beta < \kappa(\Lambda)$, the existence of two constants $c_\beta$ and $\gamma_\beta$, independent on~$\sigma$, such that:
$$ \sE_0\left[(N_0)^\beta \Ind{T_1 < \bar H_{U(x)}}\right] \pp c_\beta \ell(x)^{\gamma_\beta}. $$
where $\ell(x) = \left(\frac{\kappa}{c} \ln(x)\right)^2$ and $U(x)$ is a box centred at~$0$ with width $2\ell(x)$. We use a very similar proof to lemma~16 in \cite{perrel2024limit}, and in the following we refer to it whenever we can.

We work on the graph $\cG_\dagger = (U(x) \cup \ens{\dagger}, E_\dagger)$ obtained from $U(x)$ by collapsing all points outside of $U(x)$ in a single point $\dagger$ (therefore all the points of the boundary of $U(x)$ are linked by an edge to $\dagger$) and adding an edge from $\dagger$ to $0$.

Let $\theta^\sigma: E_\dagger \rightarrow \RR_+^{E_\dagger}$ be a function on the edges of $\cG_\dagger$ such that:
$$\Div(\theta^\sigma) = \delta_0 - \delta_\dagger. $$
We perform a change of measure from $\PP^{(\alpha)}$ to $\PP^{(\alpha + \gamma\theta^\sigma)}$. Denoting for every $\xi: E_\dagger \rightarrow \RR_+^{E_\dagger}$
\begin{align*}
Z_\xi = \prod_{u \in U(x)\cup \ens{\dagger}}\frac{\prod_{\underline{e}=u} \Gamma(\xi(e))}{\Gamma\left(\sum_{\underline{e}=u} \xi(e)\right)} \quad \text{and} \quad \omega^\xi = \prod_{e\in E_\dagger} \omega(e)^{\xi(e)},
\end{align*}
we have:
\begin{align*}
\sE_0^{(\alpha)}\left[(\omega_\sigma N_0 \Ind{T_1 < \bar H_{U(x)}})^{\beta}\right]
\pp \frac{Z_{\alpha+\gamma\theta^\sigma}}{Z_\alpha} \sE_0^{(\alpha+\gamma\theta^\sigma)}\left[(\omega_\sigma N_0)^{\beta} \omega^{-\gamma\theta^\sigma}\right]
\end{align*}

Let $a > 1$ and $b > 1$ be two conjugated exponents; by H\"older's inequality:
\begin{align*}
\sE_0^{(\alpha)}\left[(\omega_\sigma N_0 \Ind{T_1 < \bar H_{U(x)}})^{\beta}\right]
\pp \frac{Z_{\alpha+\gamma\theta^\sigma}}{Z_\alpha}  \sE_0^{(\alpha+\gamma\theta^\sigma)}\left[(N_0)^{b\beta}\right]^{1/b} \sE_0^{(\alpha+\gamma\theta^\sigma)}\left[(\omega_\sigma)^{a\beta} \omega^{-a\gamma\theta^\sigma}\right]^{1/a}.
\end{align*}
We examine one after the other the three factors of this bound.

The first expectation, by a similar argument as in lemma~16 of \cite{perrel2024limit}, is finite as soon as:
$$b\beta < \gamma$$
and can be bounded by a finite constant independent on~$\sigma$.

The second expectation can be controlled, as in lemma~16 in \cite{perrel2024limit}, by the ``max-flow min-cut theorem'': it is finite if we can find some conductances $(c^\sigma(e))_{e\in E_\dagger}$ such that:
$$ \forall e\in E_\dagger, \quad c^\sigma(e) \pp \frac{1}{\gamma(a-1)} \left[\alpha(e) + a\beta \Ind{e\in\sigma}\right], $$
and
$$ m(c^\sigma):= \inf\ens{\sum_{\partial_+A} c^\sigma(e), A \text{ connected set containing 0}} \pg 1. $$
Let us introduce a new parameter $\delta > 0$ and set:
$$ c^\sigma(e) = \frac{1}{\delta} \left[\alpha(e) + a\beta \Ind{e\in\sigma}\right]. $$
The first condition is satisfied for:
$$\delta (a-1) \pp 1.$$
As for the second condition, let $A$ be a minimizer of $m(c^\sigma)$: either $e \in \partial_+A$, and $m(c^\sigma) \pg \delta[\min \alpha + a\beta]$; or $e \not\in \partial_+A$, so $\sigma \subset A$ and $m(c^\sigma) \pg \delta\kappa(\Lambda)$. The second condition is thus satisfied if:
$$ \min\ens{\delta[\min \alpha + a\beta], \delta\kappa(\Lambda)} \pg 1.$$
If this can be achieved, the second expectation can be bounded by a finite constant independent on~$\sigma$.

We now choose the parameters $\gamma, a, b, \delta$ so as to fulfil all the conditions deduced in the two previous paragraphs:
\begin{enumerate}
\item we fix $a > 1$ such that $a\beta < \kappa(\Lambda)$---this is possible since $\beta < \kappa(\Lambda)$---and $b > 1$ its conjugated exponent;
\item we fix $\delta \in \left(\max\ens{\frac{1}{\kappa(\Lambda)}, \frac{1}{a\beta+\min\alpha}}, \frac{1}{a\beta}\right)$---which is possible thanks to the choice of~$a$;
\item we fix $\gamma \in \left(\frac{a\beta}{a-1},\frac{1}{\delta(a-1)}\right)$---which is possible thanks to the choice of~$\delta$.
\end{enumerate}

It remains to control the ratio $\frac{Z_{\alpha+\gamma\theta^\sigma}}{Z_\alpha}$. A series expansion to the order~2 of this ratio, similar to the one performed in lemma~16 of \cite{perrel2024limit}, shows that it can be bounded by some $A e^{a \|\theta^\sigma \|^2}$. By our choice of conductances, $\|\theta^\sigma \|^2$ is of order $\ln \ell(x)$, independently on~$\sigma$.

As a conclusion, there exist two constants $c_\beta$ and $\gamma_\beta$, independent on~$\sigma$, such that:
$$ \sE_0\left[(N_0)^\beta \Ind{T_1 < \bar H_{U(x)}}\right] \pp c_\beta \ell(x)^{\gamma_\beta}. $$
\end{proof}

When $\Lambda$ is big enough, $\kappa(\Lambda) > 1$. Fix $\Lambda$ such that this inequality is satisfied. It enables to have finite ${\sE_{N,0}^{(\alpha,\epsilon)}}\left[(\omega_\sigma)^{ap} N_0(H_\partial)^{ap}\right]$ for $ap < \kappa(\Lambda)$, and finite $ \EE^{(\alpha,\epsilon)}_N \left[\left(\beta_\partial \frac{\bar\pi_N^\omega(0)}{\pi_N^\omega(\partial)}\right)^a\right] $ for $a \in (1,\kappa(\Lambda))$. From this point, the reasoning about $(\QQ_N^{(\alpha,\epsilon)})_{N\pg 1}$ in the original setting is valid for $(\bar\QQ_N^{(\alpha,\epsilon)})_{N\pg 1}$ in the accelerated setting, with the difference that $\kappa$ (which can be $\pp 1$) can be replaced by $\kappa(\Lambda) > 1$: the invariant measure $\bar\QQ$ from the point of view of the accelerated particle arises as the limit in $L^1(\PP^{(\alpha)})$ of a converging subsequence of $(\bar\QQ_N^{(\alpha,\epsilon)})$ for $N \to \infty$ and $\epsilon \to 0$.

\section*{Appendix}

\begin{lemma}
\label{THappendix}
Let $(\Pi_i)_{i\pg1}$ a sequence of i.i.d. non-negative random variables, and $k$ and $N$ two integers. For every $1 \pp p \pp k$, fix $I^{(N)}_p$ a subset of $\NN$ with $N$~elements, and set $Z^{(N)}_p = \prod_{i \in I_p} \Pi_i$. Then:
$$ \EE\left[\prod_{p=1}^k Z^{(N)}_p\right] \pp \EE\left[{\Pi_1}^k\right]^N. $$

In other words, this expectation is maximal when all the $I_p$'s are the same.
\end{lemma}

\begin{proof}
First, reorganize the factors of the product of the $Z^{(N)}_p$'s in order to gather the $\Pi_i$'s:
\begin{align*}
\EE\left[\prod_{p=1}^k Z^{(N)}_p\right]
= \EE\left[\prod_{p=1}^k \prod_{i \in I^{(N)}_p} \Pi_i \right]
=\EE\left[ \prod_{i \pg 1} {\Pi_i}^{\#\ens{1 \pp p \pp k \tq i \in I_p^{(N)}}}\right]
= \prod_{i \pg 1} \EE\left[{\Pi_1}^{\#\ens{1 \pp p \pp k \tq i \in I_p^{(N)}}}\right],
\end{align*}
using the fact that the $\Pi_i$'s are i.i.d. Let us partition the index of the last product according to the value of the exponent $N_i=\#\ens{1 \pp p \pp k \tq i \in I_p^{(N)}}$ in the expectation:
$$\EE\left[\prod_{p=1}^k Z^{(N)}_p\right] = \prod_{j = 1}^k \EE\left[{\Pi_1}^j\right]^{n_j} $$
where $ n_j = \sum_{i\pg 0}\Ind{N_i=j} $ is the number of indices $i \pg 1$ belonging to exactly $j$ subsets $I_p^{(N)}$. 

By Jensen's inequality ($x \longmapsto x^{k/j}$ being convex on~$\RR_+$ when $1 \pp j \pp k$):
$$ \EE\left[\prod_{p=1}^k Z^{(N)}_p\right]
= \prod_{j = 1}^k \EE\left[{\Pi_1}^j\right]^{\frac{k}{j} \frac{j n_j}{k}} \pp \prod_{j = 1}^k \EE\left[{\Pi_1}^k\right]^{\frac{jn_j}{k}} = \EE\left[{\Pi_1}^k\right]^n $$
where:
$$  n = \sum_{j=1}^k \frac{j n_j}{k} = \frac{1}{k} \sum_{j=1}^k \sum_{i\pg 0} j \Ind{N_i=j} = \frac{1}{k} \sum_{j=1}^k \sum_{i\pg 0} N_i \Ind{N_i=j} = \frac{1}{k} \sum_{i \pg 0} N_i = N.$$

Finally $n = N$ and:
$ \EE\left[\prod_{p=1}^k Z^{(N)}_p\right]
\pp \EE\left[{\Pi_1}^k\right]^N. $
\end{proof}

\noindent
{\it\bf Acknowledgment.} The second author thanks Gady Kozma for useful discussions. Indeed, C. Sabot and G. Kozma discussed the question of the non-existence of the invariant measure 10 years ago. Although the approach was very different and unfortunately not successful, the partial answer to that question proposed here is partly motivated by these discussions. 

\noindent
{\it\bf Funding.} This work is supported by the  project  ANR LOCAL (ANR-22-CE40-0012-02)  operated by the Agence Nationale de la Recherche (ANR), and by the Institut Universitaire de France

\bibliographystyle{plain}
\bibliography{invariance_bib}

\begin{thebibliography}{10}

\bibitem{arista2024matsumoto}
Jonas Arista, Elia Bisi, and Neil O'Connell.
\newblock {Matsumoto--Yor and Dufresne type theorems for a random walk on
  positive definite matrices}.
\newblock In {\em Annales de l'Institut Henri Poincare (B) Probabilites et
  statistiques}, volume~60, pages 923--945. Institut Henri Poincaré, 2024.

\bibitem{bacallado2023edge}
Sergio Bacallado, Christophe Sabot, and Pierre Tarr\`es.
\newblock {The *-Edge-Reinforced Random Walk}.
\newblock {\em preprint, arXiv:2102.08984}, 2023.

\bibitem{berger2013local}
Noam Berger, Moran Cohen, and Ron Rosenthal.
\newblock {Local limit theorem and equivalence of dynamic and static points of
  view for certain ballistic random walks in i.i.d. environments}.
\newblock {\em The Annals of Probability}, 44(4):2889--2979, 2016.

\bibitem{berger2014effective}
Noam Berger, Alexander Drewitz, and Alejandro~F. Ramírez.
\newblock Effective polynomial ballisticity conditions for random walk in
  random environment.
\newblock {\em Communications on Pure and Applied Mathematics},
  67(12):1947--1973, 2014.

\bibitem{berman1994matrices}
Abraham Berman and Robert~J. Plemmons.
\newblock {\em {Nonnegative Matrices in the Mathematical Sciences}}.
\newblock Society for Industrial and Applied Mathematics, Philadelphia, 1994.

\bibitem{bolthausen2002static}
Erwin Bolthausen and Alain-Sol Sznitman.
\newblock {On the Satic and Dynamic Points of View for Certain Random Walks in
  Random Environment}.
\newblock {\em Methods and Applications of Analysis}, 9(3):345--376, 2002.

\bibitem{Bouchet_Sabot_Ramirez}
\'Elodie Bouchet, Alejandro~F. Ram\'irez, and Christophe Sabot.
\newblock Sharp ellipticity conditions for ballistic behavior of random walks
  in random environment.
\newblock {\em Bernoulli}, 22(2):969--994, 2016.

\bibitem{guerra2020proof}
Enrique Guerra and Alejandro~F. Ramírez.
\newblock {A proof of Sznitman's conjecture about ballistic RWRE}.
\newblock {\em Communications on Pure and Applied Mathematics},
  73(10):2087--2103, 2020.

\bibitem{komorowski2012fluctuations}
Tomasz Komorowski, Claudio Landim, and Stefano Olla.
\newblock {\em {Fluctuations in Markov processes: time symmetry and martingale
  approximation}}, volume 345.
\newblock Springer Science \& Business Media, 2012.

\bibitem{kozlov1985method}
Sergei~Mikhailovich Kozlov.
\newblock The method of averaging and walks in inhomogeneous environments.
\newblock {\em Russian Mathematical Surveys}, 40(2):73, 1985.

\bibitem{kozma2017central}
Gady Kozma and B{\'a}lint T{\'o}th.
\newblock {Central limit theorem for random walks in doubly stochastic random
  environment: ${\mathcal{H}_{-1}}$ suffices}.
\newblock {\em The Annals of Probability}, 45(6B):4307--4347, 2017.

\bibitem{Lawler82}
Gregory~F. Lawler.
\newblock Weak convergence of a random walk in a random environment.
\newblock {\em Comm. Math. Phys.}, 87(1):81--87, 1982/83.

\bibitem{bouchet2013subballistic}
Élodie Bouchet.
\newblock {Sub-ballistic random walk in Dirichlet environment}.
\newblock {\em Electronic Journal of Probability}, 18(none):1--25, 2013.

\bibitem{matsumoto2001analogue}
Hiroyuki Matsumoto and Marc Yor.
\newblock {An Analogue of Pitman's 2M-X Theorem for Exponential Wiener
  Functionals Part II: The Role of the Generalized Inverse Gaussian Laws}.
\newblock {\em Nagoya Mathematical Journal}, 162:65--86, 2001.

\bibitem{peretz2022environment}
Tal Peretz.
\newblock {Environment Viewed from the Particle and Slowdown Estimates for
  Ballistic RWRE on $\mathbb{Z}^2$ and $\mathbb{Z}^3$}.
\newblock {\em arXiv preprint arXiv:2204.12238}, 2022.

\bibitem{perrel2024limit}
Adrien Perrel.
\newblock {Limit theorem for subdiffusive random walk in Dirichlet random
  environment in dimension $ d \ge 3$}.
\newblock {\em arXiv preprint arXiv:2404.06502}, 2024.

\bibitem{poudevigne2019limit}
R\'emy Poudevigne.
\newblock {Limit Theorem for sub-ballistic random walks in Dirichlet
  environment in dimension $ d \ge 3$}.
\newblock {\em arXiv preprint arXiv:1909.03866}, 2019.

\bibitem{rapenne2023continuous}
Valentin Rapenne and Christophe Sabot.
\newblock {A continuous random operator associated with the Vertex Reinforced
  Jump Process on the circle and the real line}.
\newblock {\em arXiv preprint arXiv:2308.01120}, 2023.

\bibitem{rassoul2009almost}
Firas Rassoul-Agha and Timo Sepp{\"a}l{\"a}inen.
\newblock Almost sure functional central limit theorem for ballistic random
  walk in random environment.
\newblock In {\em Annales de l'IHP Probabilités et statistiques}, volume~45,
  pages 373--420, 2009.

\bibitem{sabot2011random}
Christophe Sabot.
\newblock Random walks in random dirichlet environment are transient in
  dimension $d \pg 3$.
\newblock {\em Probability theory and related fields}, 151(1-2):297--317, 2011.

\bibitem{sabot2013random}
Christophe Sabot.
\newblock {Random Dirichlet environment viewed from the particle in dimension
  $d\ge3$}.
\newblock {\em The Annals of Probability}, 41(2):722--743, 2013.

\bibitem{svrjp}
Christophe Sabot and Pierre Tarr{\`e}s.
\newblock {The *-Vertex-Reinforced Jump Process}.
\newblock {\em preprint, arXiv:2102.08988}, 2023.

\bibitem{sabot2011reversed}
Christophe Sabot and Laurent Tournier.
\newblock Reversed dirichlet environment and directional transience of random
  walks in dirichlet environment.
\newblock {\em Annales de l'IHP Probabilités et statistiques}, 47(1):1--8,
  2011.

\bibitem{sabot2017random}
Christophe Sabot and Laurent Tournier.
\newblock Random walks in dirichlet environment: an overview.
\newblock {\em Annales de la Faculté des sciences de Toulouse:
  Mathématiques}, 26(2):463--509, 2017.

\bibitem{Slonim2024}
Daniel~J. Slonim.
\newblock Random walks in {D}irichlet environments on {$\Bbb Z$} with bounded
  jumps.
\newblock {\em Ann. Inst. Henri Poincar\'e{} Probab. Stat.}, 60(2):1334--1355,
  2024.

\bibitem{sznitman2001class}
Alain-Sol Sznitman.
\newblock On a class of transient random walks in random environment.
\newblock {\em The Annals of Probability}, 29(2):724--765, 2001.

\bibitem{sznitman2002effective}
Alain-Sol Sznitman.
\newblock An effective criterion for ballistic behavior of random walks in
  random environment.
\newblock {\em Probability theory and related fields}, 122(4):509--544, 2002.

\bibitem{Sznitman_ten_lectures}
Alain-Sol Sznitman.
\newblock {\em Lectures on random motions in random media { In } Ten lectures
  on random media}, volume~32 of {\em DMV Seminar}.
\newblock Birkh\"auser Verlag, Basel, 2002.

\bibitem{sznitman2003new}
Alain-Sol Sznitman.
\newblock On new examples of ballistic random walks in random environment.
\newblock {\em The Annals of Probability}, 31(1):285--322, 2003.

\bibitem{sznitman1999law}
Alain-Sol Sznitman and Martin Zerner.
\newblock A law of large numbers for random walks in random environment.
\newblock {\em The Annals of Probability}, 27(4):1851--1869, 1999.

\bibitem{toth2018quenched}
B{\'a}lint T{\'o}th.
\newblock Quenched central limit theorem for random walks in doubly stochastic
  random environment.
\newblock {\em The Annals of Probability}, 46(6):3558--3577, 2018.

\bibitem{tournier2009integrability}
Laurent Tournier.
\newblock {Integrability of exit times and ballisticity for random walks in
  Dirichlet environment}.
\newblock {\em Electronic Journal of Probability}, 14(none):431--451, 2009.

\bibitem{tournier2015asymptotic}
Laurent Tournier.
\newblock Asymptotic direction of random walks in dirichlet environment.
\newblock In {\em Annales de l'IHP Probabilités et statistiques}, volume~51,
  pages 716--726, 2015.

\bibitem{yilmaz2010averaged}
Atilla Yilmaz.
\newblock Averaged large deviations for random walk in a random environment.
\newblock {\em Annales de l'IHP Probabilit{\'e}s et statistiques},
  46(3):853--868, 2010.

\bibitem{zeitouni2004random}
Ofer Zeitouni.
\newblock Random walks in random environment.
\newblock {\em Lecture notes in Mathematics}, 1837:190--312, 2004.

\end{thebibliography}
	
\end{document}